\documentclass[	a4paper, 
								11pt]{article}

\usepackage{amssymb,amsmath,amsthm}
\usepackage{graphicx}
\usepackage{subfigure}
\usepackage{overpic}
\usepackage[table]{xcolor}

\usepackage{chngcntr}
\usepackage{lscape}
\usepackage{rotating}
\usepackage{booktabs}
\usepackage{units}

\usepackage{epstopdf}

\usepackage{multirow}
\usepackage{calc}

\usepackage{tikz}
\usepackage{tikz-cd}
\usepackage{pgfplotstable}
\usepackage{multicol}
\usepackage{pgfplots}
\usetikzlibrary{3d,calc}
\usetikzlibrary{shapes, spy, scopes, matrix, datavisualization.formats.functions, arrows, datavisualization, plothandlers, plotmarks, svg.path, positioning}
\usetikzlibrary{external}
\tikzset{external/only named=true}
\pgfplotsset{compat=newest}
\pgfplotsset{plot coordinates/math parser=false}
\pgfplotsset{scaled y ticks = false, tick label style={/pgf/number format/fixed}}
\pgfplotsset{scaled x ticks = false, tick label style={/pgf/number format/fixed}}

\usepackage{float}

\usepackage[bookmarksopen, pdfstartview={FitV}, linkbordercolor=white]{hyperref}

\usepackage{nomencl}

\usepackage[]{algorithm2e}
\usepackage{todonotes}

\usepackage[capitalize,nameinlink]{cleveref}
	
\newtheoremstyle{mystyle}
  {}
  {}
  {\itshape}
  {}
  {\bfseries}
  {.}
  { }
  {}%
	
\theoremstyle{mystyle}
\newtheorem{theorem}{Theorem}[section]
\theoremstyle{plain}

\theoremstyle{plain}

\theoremstyle{plain}

\theoremstyle{plain}

\theoremstyle{plain}
\newtheorem{corollary}[theorem]{Corollary}
\theoremstyle{plain}
\newtheorem{proposition}[theorem]{Proposition}

\newcommand{\BIGOP}[1]{\mathop{\mathchoice%
{\raise-0.22em\hbox{\huge $#1$}}%
{\raise-0.05em\hbox{\Large $#1$}}{\hbox{\large $#1$}}{#1}}}

\newcommand{\BIGboxplus}{\mathop{\mathchoice%
{\raise-0.35em\hbox{\huge $\boxplus$}}%
{\raise-0.15em\hbox{\Large $\boxplus$}}{\hbox{\large $\boxplus$}}{\boxplus}}}

\makeatletter
\def\tagform@#1{\maketag@@@{\ignorespaces#1\unskip\@@italiccorr}}
\let\orgtheequation\theequation
\def\theequation{(\orgtheequation)}
\makeatother

\def\N{\mathbb N}			 \def\R{\mathbb R}

\def\R{\mathbb{R}}
\def\N{\mathbb{N}}

\def\cH{\mathcal{H}}

\def\cK{\mathcal{K}}

\def\cO{\mathcal{O}}

\def\ee{\boldsymbol{e}}
\def\ff{\boldsymbol{f}}

\def\kk{\boldsymbol{k}}

\def\xx{\boldsymbol{x}}
\def\yy{\boldsymbol{y}}
\def\zz{\boldsymbol{z}}

\def\AA{\boldsymbol{A}}
\def\BB{\boldsymbol{B}}
\def\CC{\boldsymbol{C}}

\def\II{\boldsymbol{I}}

\def\LL{\boldsymbol{L}}
\def\MM{\boldsymbol{M}}

\def\QQ{\boldsymbol{Q}}

\def\TT{\boldsymbol{T}}

\def\II{\boldsymbol{I}}
\def\QQ{\boldsymbol{Q}}

\def\11{\mathbf{1}}

\def\00{\boldsymbol{0}}

\def\aalpha{\boldsymbol{\alpha}}
\def\bbeta{\boldsymbol{\beta}}

\def\oomega{\boldsymbol{\omega}}

\def\transpose{\text{T}}
\def\spn{\operatorname{span}}
\def\d{\operatorname{d}\!}

\usepackage[left=2.5cm,right=2.5cm,top=2.5cm,bottom=2.5cm]{geometry}
\theoremstyle{plain}

\usepackage[toc]{appendix}
\usepackage[export]{adjustbox}
\newcommand{\diam}{\operatorname{diam}}

\usepackage{authblk}

\usepackage{comment}

\def\Letters{A,B,C,D,E,F,G,H,I,J,K,L,M,N,O,P,Q,R,S,T,U,V,W,X,Y,Z}
\makeatletter
\@for \@l:=\Letters \do{%
  \expandafter\edef\csname\@l bb\endcsname{%
  \noexpand\ensuremath{\noexpand\mathbb{\@l}}}%
  \expandafter\edef\csname\@l bf\endcsname{{\noexpand\bf \@l}}%
  \expandafter\edef\csname\@l cal\endcsname{%
  \noexpand\ensuremath{\noexpand\mathcal{\@l}}}%
  \expandafter\edef\csname\@l eu\endcsname{%
  \noexpand\ensuremath{\noexpand\EuScript{\@l}}}%
  \expandafter\edef\csname\@l frak\endcsname{%
  \noexpand\ensuremath{\noexpand\mathfrak{\@l}}}%
  \expandafter\edef\csname\@l rm\endcsname{{\noexpand\rm \@l}}%
  \expandafter\edef\csname\@l scr\endcsname{%
  \noexpand\ensuremath{\noexpand\mathscr{\@l}}}%
}
\newcommand{\bs}[1]{{\boldsymbol#1}}
\newcommand{\isdef}{\mathrel{\mathrel{\mathop:}=}}
\newcommand{\defis}{\mathrel{=\mathrel{\mathop:}}}

\title{On Quasi-Localized Dual Pairs in Reproducing Kernel Hilbert Spaces}
\date{\today}

\author[1]{Helmut Harbrecht}
\author[2]{R\"udiger Kempf}
\author[3]{Michael Multerer}
\affil[1]{
	   Departement Mathematik und Informatik,
	   Universit\"at Basel,
	   Spiegelgasse 1, 
      4051 Basel,
	   Switzerland}
\affil[2]{Applied and Numerical Analysis,
 	   Department of Mathematics,
 	   University of Bayreuth, 
 	   95440 Bayreuth,
 	   Germany}   
\affil[3]{
	   Istituto Eulero,
	   Universit\`a della Svizzera italiana,
	   Via la Santa 1, 
      6962 Lugano,
	   Switzerland}
\begin{document}
\maketitle

\abstract{In scattered data approximation, the span of a finite
number of translates of a chosen radial basis function is used 
as approximation space and the basis of translates is
used for representing the approximate. However, this natural
choice is by no means mandatory and different choices, like, 
for example, the Lagrange basis, are possible and 
might offer additional features. In this article, we discuss 
different alternatives together with their canonical duals.
We study a localized version of the
Lagrange basis, localized orthogonal bases, such as the Newton basis,
and multiresolution versions thereof, constructed by means
of samplets. We argue that the choice of orthogonal bases
is particularly useful as they lead to symmetric preconditioners.
All bases under consideration are compared numerically to 
illustrate their feasibility for scattered data approximation.
We provide benchmark experiments in two spatial dimensions and 
consider the reconstruction of an implicit surface as
a relevant application from computer graphics.}

\section{Introduction}
Kernel methods have become more and more popular over the years. Applications 
range from scattered data approximation to machine learning, see for example
\cite{Fasshauer2007,Christmann:SVM,wendland:ScatteredDataApproximation} and 
the references therein. In general, these methods aim at the reconstruction 
of an unknown function by only using scattered data, i.e., tuples of data sites 
and corresponding measurements of the function. Assuming that the data generating
process $ f $ is contained in a reproducing kernel Hilbert space $ \cH $, we can 
always find a set of \emph{dual pairs} $\big\{\big( \phi_i, \widetilde{\phi}_i\big)\big\} $ 
such that the orthogonal projection $s_f$ of $f$ in $ \operatorname{span}\{ \phi_i \}$ can 
be written as 
\begin{align*}
   s_f=\sum_{i=1}^{N}\big\langle f, \widetilde{\phi}_i\big\rangle_{\cH} \phi_i \quad \text{or, equivalently,} \quad 
   s_f= \sum_{i=1}^{N} \langle f,  \phi_i \rangle_{\cH} \widetilde{\phi}_i. 
\end{align*}
Given a set of data sites \(\{\xx_1,\ldots,\xx_N\}\),
the most commonly known dual pair is $ \{(K(\cdot, \xx_i), \chi_i)\} $. Herein,
$\{K(\cdot, \xx_i)\}$, is the basis of kernel translates induced by the reproducing kernel
$K$ of $ \cH $ and $ \chi_i $ is the corresponding
Lagrange basis. This specific pair allows the easy representation of the orthogonal projection by the
interpolant
\begin{align*}
    s_f=\sum_{i=1}^N f(\xx_i) \chi_i.
\end{align*}
This representation has the advantage that it separates the data from the data sites. 
However, even if we do not consider the cost to compute the Lagrange basis
$\{\chi_i\}$, every evaluation of the interpolant \(s_f\) entails a cost of $ \cO(N^2) $, since every 
basis element
$ \chi_i $ is a linear combination of the $ N $ kernel translates that need to be evaluated. 
Numerically, this is not feasible. Therefore, a framework was established in 
a series of articles \cite{hangelbroek:InverseTheoremLocalizedLagrangeFunctions,
fuselier:LocalizedLagrangeFunctions,hangelbroek:LpStabilityLagrangeBasis}
in order to compute localized approximate Lagrange functions for very
specific classes of kernels. 

In the present article, we give an interpretation of the localization 
of the Lagrange basis by considering it as an approximation/sparsification 
of the inverse of the kernel matrix \({\bs A}=[K({\xx_i,\xx_j})]\)
by means of a very specific pattern, 
which we call \emph{footprint}. This approximate inverse is typically not 
symmetric although the kernel matrix itself is. Hence, we cannot use it 
as a preconditioner in a conjugate gradient method. To address this issue, 
we discuss modifications of the approach by involving the matrix
square root of the footprints or their Cholesky decomposition. In the limit, for growing
footprint size, both approaches yield orthonormal basis functions, which means that 
$\tilde\phi_i = \phi_i$ for $i=1,\ldots,N$. Particularly, the second approach
is a localized version of the Newton basis, see \cite{mueller:NewtonBasisForKernelSpaces}.
Either choice may then be used to construct symmetric and sparse preconditioners 
for notoriously bad conditioned kernel matrices, e.g., issuing from the widely 
used Mat\'ern kernels.

A multiresolution approach for the construction of dual pairs, that we consider here,
are samplets \cite{HM2}. Samplets are discrete signed measures
constructed such that polynomials up to a degree of our choice vanish. Therefore, the 
kernel matrix becomes quasi-sparse in samplet coordinates and can be compressed to
a sparse matrix. The resulting matrix pattern is finger band structured with easy
to compute Cholesky decomposition whence a graph based reordering strategy is used 
\cite{HM2}. Furthermore, it has been shown in \cite{HMSS} that samplets establish 
even a sparse arithmetic for kernel matrices, enabling efficient scattered data 
analysis.

This article is organized as follows. In \cref{sec:ScatteredDataApproximation}, 
we recall the basics of scattered data approximation and the setting 
we are studying. In \cref{sec:Lagrange}, we recall general dual pairs and
the Lagrange basis in particular. Moreover, we discuss their connection to 
the theory of pseudo-differential operators. This gives an alternative point 
of view to important properties of kernel matrices which we want to use 
in the subsequent section. There, in \cref{sec:LocLagrangeFunctions}, we 
consider alternative dual pairs. First, in \cref{subsec:LocLagrangeFunctions}, 
we give a review of the quasi-localization of the Lagrange basis by using 
the exponential decay of the inverse kernel matrix. Then, we use these ideas 
to derive a new way to construct preconditioners for kernel matrices in 
\cref{subsec:LocLagrangeFunctions2}. The underlying basis functions are 
orthonormal, including especially a localized version of the Newton basis.
As an alternative approach, we review the basics of samplets in \cref{sec:Samplets}. 
In \cref{sec:Numerics}, we compare the different dual pairs under consideration
numerically. Extensive numerical tests in two and three spatial dimensions
are provided. Finally, concluding remarks are stated in \cref{sec:conclusio}.


\section{Scattered Data Approximation}\label{sec:ScatteredDataApproximation}

Let $ \Omega \subseteq \R^d $ denote a domain, and 
\(\Hcal\) Hilbert space of functions mapping \(\Omega\) to \(\Rbb\).
If $ \cH \subseteq C(\Omega) $, then \(\Hcal\) is a \emph{reproducing kernel Hilbert 
space} (RKHS) and there exists a unique \emph{reproducing kernel} 
$ K\colon\Omega \times \Omega \to \R $ such that $ K(\cdot, \xx) \in \cH $ for all 
$ \xx \in \Omega $ and each element $ f \in \cH $ can be point-wise 
recovered, i.e., $ f(\xx) = \langle f, K(\cdot, \xx) \rangle_{\cH}$. 
It is well-known, see, e.g.\ \cite{wendland:ScatteredDataApproximation}, 
that reproducing kernels are symmetric and \emph{positive semi-definite}, 
i.e., the kernel matrix $ \AA = [K(\xx_i,\xx_j)]_{i,j}$ is symmetric and 
positive semi-definite for all $ \{ \xx_1, \dots, \xx_N \} \subseteq \Omega $ 
and all $ N \in \N $.

In this article, we focus on kernels that are induced by
\emph{radial functions}. A function \(\Phi\colon\Rbb^d\to\Rbb\) is called
radial, iff there is a function $ \phi: [0, \infty) \to \R $ such that \(\Phi(\xx)=\phi(\|\xx\|_2)\).
If the function \(\Phi\) is positive definite in the sense that the
induced kernel \(K(\xx,\yy)\isdef\Phi(\xx-\yy)\) is positive definite, then
it gives rise to a reproducing kernel.

One of the most popular of these radial functions is known as 
\emph{Mat\'ern kernel} or \emph{Sobolev spline} 
$\Phi_{\nu} : \R^d \to \R $, dependent on the \emph{smoothness parameter}
$ \nu > d/2 $. 
It is defined by 
\begin{align}
 \phi_{\nu}(r) = \frac{2^{1-\nu}}{\Gamma(\nu)} r^{\nu - \frac{d}{2}} K_{\nu - \frac{d}{2}} (r),
 \quad  r \geq 0,
\end{align}
where $ \Gamma $ is the Riemannian Gamma function and $ K_{\beta} $ 
is the modified Bessel function of the second kind, see \cite{MAT} for example.
For specific values of the smoothness parameter
$ \nu $ the representation of $ \phi_{\nu} $ simplifies 
significantly, we give a selection of specific representations in \cref{tab:exampleMatern}.

\begin{table}[hbt]
\begin{center}
   \begin{tabular}{|c|c|c|}\hline
    $ \nu $ & $ \phi_{\nu} (r) $ & smoothness \\
    \hline
    $ 1/2 $ & $ \exp(-r)  $ & $ C^0 $\\
    $ 3/2 $ & $ (1+r)\exp(-r) $ & $ C^2 $\\
    $ 5/2 $ & $ (3+3r+r^2)\exp(-r) $ & $ C^4 $ \\
    $ \infty $ & $ \exp(-r^2) $ & $ C^{\infty} $ \\\hline
\end{tabular} 
\caption{Examples for different smoothness parameters 
$\nu$ of the Mat{\'e}rn kernel.}
\label{tab:exampleMatern}
\end{center}
\end{table}

It is well-known that the Fourier transform of $ \Phi_{\nu} $ 
decays algebraically, i.e., there holds
\begin{align*}
\widehat{\Phi_{\nu}}(\oomega) = (1 + \| \oomega \|_2^2)^{- \nu - \frac{d}{2}}, \quad \oomega \in \R^d,
\end{align*}
see \cite[Theorem 6.13]{wendland:ScatteredDataApproximation} for example.
This means in particular, see, e.g., \cite{wendland:ScatteredDataApproximation}, 
that $ \Phi_{\nu} $ is the reproducing kernel of the Sobolev-Hilbert space 
$ H^{\nu - \frac{d}{2}}(\R^d) $. Especially, the kernel defines a 
pseudo-differential operator on $\R^d$, a fact which we will 
exploit later on.

Having fixed the kernel of interest, we now describe 
the problem of function approximation. We are interested 
in recovering an unknown function $ f \in \cH $, where 
we are given only a finite set of data $ \{(\xx_N, f_N),\ldots,(\xx_N, f_N)\}
\subset \Omega \times \R$. We collect 
the abscissae $ \xx_i $ in a \emph{set of data sites} $ X\isdef \{ \xx_1, 
\dots, \xx_N \} \subset \Omega $. Associated to this 
set are two characteristic quantities: the \emph{fill distance} 
\begin{align*}
h_{X,\Omega} := \sup_{\xx \in \Omega} \min_{\xx_i \in X} \| \xx - \xx_i \|_2
\end{align*}
and the \emph{separation radius}
\begin{align*}
q_X := \frac{1}{2} \min_{i \neq j} \| \xx_i - \xx_j \|_2.
\end{align*}
For the theoretical results we present later, we require that the set of 
data sites is \emph{quasi-uniform}, i.e., 
that there is a constant $ c_{\operatorname{qu}} > 0 $ such that 
\begin{align*}
q_X \leq h_{X,\Omega} \leq c_{\operatorname{qu}} q_X.
\end{align*}
An easy comparison of volumes then yields, that there are constants $ c_1 $ and $ c_2 $ such that
\begin{align*}
    c_1 N^{- \frac{1}{d}} \leq h_{X,\Omega} \leq c_2 N^{- \frac{1}{d}},
\end{align*}
see, e.g., \cite[Proposition 14.1]{wendland:ScatteredDataApproximation}. Clearly, a similar estimate holds for the separation radius, too.

To recover the unknown function, we look at two approximation methods, \emph{interpolation} and, more generally,
\emph{regularized least squares approximation} or \emph{kernel ridge regression}.
Both can be seen as finding the minimizer of the functional $ J_{\lambda} $ within the subspace
\[
\cH_X\isdef\operatorname{span}\{K(\cdot,\xx_1),\ldots,K(\cdot,\xx_N)\}
\]
spanned by the \emph{basis of kernel translates}  within an RKHS $ \cH $.
\begin{align*}
J_{\lambda}(s)\isdef\sum_{i=1}^{N} | f_i - s(\xx_i) |^2 + \lambda \| s \|_{\cH}^2, \quad s \in \cH_X,
\end{align*}
where $ \lambda \geq 0 $ is a \emph{regularization parameter}. 
If we set $ \lambda = 0 $ we force interpolation, 
i.e., the minimizer $s_{f}\isdef s_0^*$ of $ J_{0} $ satisfies 
$ s_{f}(\xx_i) = f_i $, for $i=1\,\ldots,N$. This should 
only be chosen if we deal with no noise on the data, that is, if it holds 
$ f_i = f(\xx_i) $. In this case, the minimizer \(s_{f}\) coincides with
the orthogonal projection of \(f\) onto \(\Hcal_X\), since there holds
\(\langle f-s_{f},v\rangle_\cH=0\) for any \(v\in\cH_X\) by the reproducing property.

In any other case, one should use 
$ \lambda > 0 $. Even so, the numerical procedure 
is the same, stated in the next theorem, which is a 
version of the \emph{representer theorem}, see, 
e.g., \cite{Christmann:SVM}.

\begin{theorem}
Let $ \Omega \subseteq \R^d $ be a domain. Let 
$ K $ be the reproducing kernel of an RKHS
$ \cH $ on $ \Omega $ and $ X = \{ \xx_1, \dots, \xx_N\} \subseteq \Omega $
be a set of sites. Then, for all $ \lambda \geq 0 $, there exists 
a unique solution $ s^*_{\lambda} $ of $ \min J_{\lambda} (s) $. 
In addition, there exists a coefficient vector 
$ \aalpha \in \R^N $ such that 
\begin{align*}
s^*_{\lambda} = \sum_{j=1}^{N} \alpha_j K(\cdot , \xx_j),
\end{align*}
i.e., $ s^*_{\lambda} \in\cH_X$. Furthermore, 
$ \aalpha $ is the unique solution of the linear system
\begin{equation}\label{eq:LGS}
( \AA + \lambda \II) \aalpha = \ff,
\end{equation}
where $ \AA = [K(\xx_i, \xx_j)]_{i,j} \in \R^{N \times N} $ is
the kernel matrix, 
$ \II \in \R^{N \times N} $ is the identity matrix, and 
$ \ff = [f_i]\in \R^N $ is the data vector.
\end{theorem}

\section{Other Bases for $ \cH_X $}\label{sec:Lagrange}
By definition of the approximation space, $ \cH_X $ is the linear 
span of the basis of kernel translates. 
However, this might not be the most practical basis. For example, 
it is well-known that there is a \emph{Lagrange basis} 
$ \{ \chi_j \} $ satisfying $ \chi_j(\xx_i) = \delta_{i,j} $, $i,j=1\ldots, N$. 
The Lagrange basis is well-understood and we refer to the vast amount
of available literature, see
\cite{hangelbroek:InverseTheoremLocalizedLagrangeFunctions,
fuselier:LocalizedLagrangeFunctions,wendland:LagrangeFunctions,
wendland:ScatteredDataApproximation} and the references therein. 
Another basis is the \emph{Newton basis}, which amounts to an
orthogonal basis of \(\Hcal_X\), 
compare \cite{mueller:NewtonBasisForKernelSpaces,Pazouki2011}.
The general concept of \emph{dual pairs} is introduced in 
\cref{subsec:DualPairs}.
We repeat some important properties in \cref{subsec:FullLagrangeFunctions} 
and give an alternative proof of 
the exponential decay property by interpreting the underlying operators as 
\emph{pseudo-differential operators}. This view, although well-established 
in the theory of partial differential equation, is new in the context of 
approximation theory. 

\subsection{Dual Pairs}
\label{subsec:DualPairs}
A \emph{dual pair} $\big\{\big( \phi_i, \widetilde{\phi}_i\big)\big\}$ is a 
set of pairs of basis functions for $\cH_X$ such
that $\big\langle\phi_i,\widetilde\phi_j\big\rangle_\Hcal = \delta_{i,j}$
for all $i,j=1,\ldots,N $. Thus, we immediately obtain the representations
\begin{align*}
 s_f=   \sum_{i=1}^{N}\big\langle f, \widetilde{\phi}_i\big\rangle_{\cH} \phi_i 
    = \sum_{i=1}^{N} \langle f,  \phi_i \rangle_{\cH} \widetilde{\phi}_i
\end{align*}
for the sought kernel interpolant. Its computation can therefore be
accelerated if the basis functions under consideration admit 
sparse representations.

Obviously, the kernel matrix
\begin{align*}
  \AA = [K(\xx_i,\xx_j)]_{i,j=1}^N 
    = [\langle K(\cdot,\xx_i),
        K(\cdot,\xx_j\rangle_{\cH}]_{i,j=1}^N
\end{align*}
is the Gramian of basis of kernel translates. Hence, any choice of $ \gamma \in \R $ yields 
a dual pair by defining
\begin{align*}
\big[
\widetilde{\phi}_1,\ldots,
\widetilde{\phi}_N\big]
={\bs k}(\cdot)\AA^{-1-\gamma}
\quad\text{and}\quad
[{\phi}_1,\ldots,
{\phi}_N
]
={\bs k}(\cdot)\AA^{\gamma},
\end{align*}
where we set
\begin{equation*}
{\bs k}({\bs x})=[K(\cdot,\xx_1),\ldots,K(\cdot,\xx_N)].
\end{equation*}
Canonical choices are $ \gamma=0 $, yielding the \emph{Lagrange basis},
and $ \gamma=-1/2 $, yielding an orthogonal basis, where we have 
$\phi_i = \widetilde{\phi}_i$ for all $i=1,\ldots, N$. 
Furthermore, we notice that we have 
\begin{align*}
\AA^{-1} = \AA^{-1/2} \AA^{-1/2}=\LL \LL^{\transpose},
\end{align*}
where $ \LL $ is the Cholesky factor of $ \AA ^{-1} $. Hence, 
there exists an isometry such that $ \LL = \QQ \AA^{-1/2} $,
which implies that another orthonormal basis of the space $\cH_X$ is obtained 
by the Cholesky decomposition of $ \AA^{-1} $. It is called \emph{Newton 
basis} in literature, since the $i$-th basis function has $i$ zeros, compare 
\cite{mueller:NewtonBasisForKernelSpaces,Pazouki2011}.

\subsection{About the Lagrange Basis}
\label{subsec:FullLagrangeFunctions}
We shall first collect some results for the Lagrange basis.
Since the basis elements \(\chi_i\) are elements of \(\Hcal_X\), 
we can express them as 
 \begin{align*}
 \chi_i = \sum_{k=1}^{N} \alpha^{(i)}_k K(\cdot, \xx_k), \quad i=1,\ldots,N,
 \end{align*}
for certain coefficient vectors $ \aalpha^{(i)}\in\R^N$ 
for any $i=1,\ldots,N$. There are different ways to compute these coefficient vectors.
We consider two of them.

The first one allows us to compute these vectors a-priori in an 
offline phase, as soon as the set of data sites is fixed. The Lagrange 
condition $ \chi_j(\xx_i) = \delta_{i,j} $ directly leads to the 
linear system
\begin{align*}
\AA \aalpha^{(i)} = \ee_i
\end{align*}
for the $ i $-th coefficient vector, where $\AA  = [K(\xx_i, \xx_j)]_{i,j}$ 
is again the kernel matrix and $ \ee_i \in \R^N $ is the $ i $-th unit vector.

With the second method we can compute all $ \chi_i(\xx) $, i.e., the Lagrange basis
evaluated in a $ \xx \in \R^d $, simultaneously, by solving the linear system
\[
[\chi_1 (\xx),\ldots,\chi_N (\xx)]\AA
= \kk (\xx),
\]
which has a unique solution for every $ \xx \in \R^d $. This leads 
to the representation
\begin{align}
\chi_i(\xx) &= \sum_{j=1}^{N} [\AA^{-1}]_{i,j} K(\xx, \xx_j) \label{eq:fullLagrangeFunction1} \\
&= \kk (\xx)\AA^{-1} \ee_i , \quad i=1,\ldots, N. \label{eq:fullLagrangeFunction2}
\end{align}

With the Lagrange functions at hand, we can represent the minimizer $ s_f $ 
of $ J_0 $ by
\begin{align}\label{eq:representationInterpolantLagrangeFunctions}
s_f(\xx) = \sum_{i=1}^{N} f(\xx_i) \chi_i(\xx).
\end{align}
Often, this representation is associated with a \emph{quasi-interpolation process}, see, 
e.g., \cite{Buhmann:Quasi-Interpolation}, however, in the context of this paper, we do not 
expect to have polynomial reproduction, a property commonly expected of quasi-interpolation.

Motivated by the desire to obtain a similar representation for $ s^*_{\lambda} $ 
with $ \lambda > 0 $, we define the \emph{modified Lagrange basis} via
\begin{align}
{\chi}_{\lambda,i} (\xx) &= \sum_{j=1}^{N} \big[\AA + \lambda \II\big]^{-1}_{i,j} 
K(\xx, \xx_j) \label{eq:modifiedLagrangeFunction1} \\
&= \kk (\xx)(\AA + \lambda \II)^{-1} \ee_i , 
\quad i=1,\ldots,N.\label{eq:modifiedLagrangeFunction2}
\end{align} 
Employing the modified Lagrange basis, we can represent $ s^*_{\lambda} $ as 
\begin{align}\label{eq:representationApproximationLagrangeFunctions}
s^*_{\lambda} (\xx) = \sum_{i=1}^{N} f_i{\chi}_{\lambda,i} (\xx).
\end{align}

For the Mat\'ern kernel $ K(\xx,\yy) = \Phi_{\nu}(\xx-\yy) $, it turns out that the 
dual basis associated to the basis of kernel translates, i.e., the Lagrange basis, 
is {highly localized}, that is, $ |\chi_i (\xx)| $ decays 
exponentially fast if $ \| \xx - \xx_i \|_2 $ grows. This is a direct 
consequence of the fact that the entries of the inverse of the kernel 
matrix $ \AA $ exhibit an exponential off-diagonal decay. These can 
be quantified as follows.

\begin{corollary}\label{coro:off-diagonal decay}
Let $ X = \{ \xx_1, \dots, \xx_N \} $ be a quasi-uniform 
set of sites with fill distance $ h_{X,\Omega}  < h_0 $ and 
separation radius $ q_X $. Let $ \Phi_{\nu} $ be the Mat{\'e}rn 
kernel with smoothness parameter $ \nu $, such that $ \nu - d/2 \in \N $. Then there is a 
constant $ C > 0 $ such that the entries of the 
inverse of the kernel matrix satisfy the estimate
\begin{align}\label{eq:exponentialDecayMatrix}
\big| \big[\AA^{-1}\big]_{j,k}\big| \leq C q_X^{d - 2 \nu} 
\exp \left( - \eta \frac{\| \xx_j - \xx_k \|_2}{h_{X,\Omega}} \right),
\end{align}
where $ \eta = \eta(\nu,d) < 1 $ is a positive constant. 
This implies the following decay condition on the Lagrange basis
\begin{align*}
| \chi_i(\xx) | \leq C \left( \frac{h_{X,\Omega}}{q_X} \right)^{\nu - \frac{d}{2}} 
\exp \left( - \mu \frac{\| \xx - \xx_i \|_2}{h_{X,\Omega}} \right), \quad \xx \in \R^d,\text{ for }i=1,\ldots,N,
\end{align*}
with $ \mu = - \frac{h_0}{4} \log(\eta) $. 
\end{corollary}

The first proofs  of this corollary were given in 
\cite{fuselier:LocalizedLagrangeFunctions} for 
$ \Omega = \mathbb{S}^{d-1} $, the $ d-1 $ dimensional sphere. In 
\cite{hangelbroek:InverseTheoremLocalizedLagrangeFunctions}, 
it was refined to the version given here for general bounded 
domains. In the next subsection, we offer an alternative proof 
by considering the kernel in the context of
pseudo-differential operators.

We see in \eqref{eq:fullLagrangeFunction2} and \eqref{eq:modifiedLagrangeFunction2} 
that the main computational task to evaluate the (modified) Lagrange basis
$\{\chi_i\}$ and $\{\chi_{\lambda,i}\} $, respectively, is the inversion of the 
dense, positive definite matrix $\AA $ or $\AA + \lambda\II $, respectively,
whose entries decay exponentially apart from the diagonal. Hence, the desired fast evaluation 
of the approximation $ s^*_{\lambda} $ can be achieved by finding efficient 
ways to approximate the kernel matrix $ \AA $ and its inverse $ \AA^{-1} $.

\subsection{Alternative Point of View: 
Pseudo-differential Operators}\label{subsec:pseudos}
In this subsection, let $\Omega\subset\mathbb{R}^d$ be a smooth domain. 
Then, the 
reproducing kernel under consideration gives rise to a compact 
integral operator on \(L^2(\Omega)\) according to
\begin{equation}\label{eq:pseudo}
  \cK:L^2(\Omega)\to L^2(\Omega), 
  \quad u \mapsto \int_\Omega K(\cdot,\yy) u(\yy)\d\yy.
\end{equation}
For many kernel functions, $\cK$ constitutes a classical 
pseudo-differential operator of negative order $s < 0$. We 
refer to e.g.~\cite{HorIII,Rodino92,Seeley,Taylor81} for the 
details of this theory, including the subsequent developments.

We are especially interested in pseudo-differential operators 
$\mathcal{K}$ which belong to the subclass $OPS_{cl,1}^s$ of 
\emph{analytic} pseudo-differential operators, see \cite{Rodino92}. 
Their kernels are known to be \emph{asymptotically smooth}, 
meaning that there holds the decay property
\begin{equation}\label{eq:decay}
\begin{aligned}
&\forall\boldsymbol\alpha,\boldsymbol\beta\in\mathbb{N}^d:\ 
\exists\,\rho>0:\\
&\qquad|\partial^{\boldsymbol\alpha}_{\xx} 
\partial^{\boldsymbol\beta}_{\yy} K(\xx,\yy)|
\leq C \frac{(|\boldsymbol\alpha|+|\boldsymbol\beta|)!}
{\big(\rho
\|\xx-\yy\|_2^{s+d}\big)^{|\boldsymbol\alpha|+|\boldsymbol\beta|}},
\quad (\xx,\yy)\in(\Omega\times\Omega)\setminus\Delta
\end{aligned}
\end{equation}
apart from the diagonal $\Delta\isdef\{(\xx,\yy)\in\Omega
\times\Omega:\xx = \yy\}$ for some constant \(C>0\). 
This estimate is the key for 
the computational treatment of the underlying kernel 
matrices by means of hierarchical matrix techniques
\cite{hackbusch:HierarchicalMatrices} or wavelet matrix 
compression \cite{BCR,DHS}.

A special feature of pseudo-differential operators is
that they form an algebra, which we can be exploited 
for kernel matrices. Besides addition, the concatenation of two 
pseudo-differential operators $\mathcal{K}$ and $\mathcal{K}'$
essentially corresponds to the product of the respective kernel 
matrices $\AA$ and $\BB$ under the assumption that the set of 
data sites is asymptotically uniformly distributed modulo one.
This means, we have for every Riemann integrable function 
$f\colon\Omega\to\mathbb{R}$ that
\[
  \bigg|\int_\Omega f(\xx)\d\xx
  -\frac{|\Omega|}{N}\sum_{i=1}^N f(\xx_i)\bigg|
  \to 0\quad\text{as $N\to\infty$},
\]
compare \cite{LP10}. Under this assumption, the concatenation 
$\mathcal{K}\circ\mathcal{K}'$ corresponds to the matrix-matrix 
product $\AA\BB$. Indeed, the kernel $K''$ of the concatenated 
pseudo-differential operator satisfies
\[
  K''(\xx_i,\xx_j) = \int_\Omega K(\xx_i,\zz) K'(\zz,\xx_j)\d\zz
  \approx \frac{|\Omega|}{N}\sum_{k=1}^N K(\xx_i,\xx_k) K'(\xx_k,\xx_j)
  = [\AA]_{i,:} [\BB]_{:,j}.
\]
This means that the matrix $\AA\BB$ approximates the kernel matrix 
induced by $K''$. As this is a compressible pseudo-differential
operator, $\AA\BB$ is also compressible, compare \cite[Theorem 4]{HMSS}. 
This property can be used to show that the inverse kernel matrix is 
compressible.

Consider a symmetric and positive definite kernel 
$K$ such that the associated pseudo-differential 
operator satisfies $\cK\in OPS_{cl,1}^s$ with $s<0$. Then, the 
inverse of $\cK+\lambda I$ is of the form $\lambda^{-1} I-\cK'$ 
with $\cK'$ being likewise a pseudo-differential operator from 
$OPS_{cl,1}^s$, see \cite[Lemma 2]{HMSS}. Especially, the 
Schwarz kernel $K'$ which underlies the operator $\cK'$ by the 
Schwartz kernel theorem, see, e.g., \cite[Chap.~5]{HorI},
is symmetric, positive definite, and likewise asymptotically 
smooth. Thus, the inverse of a (regularized) kernel matrix 
is compressible in view of the above mentioned link between
the matrix-matrix product and the concatenation of the associated
pseudo-differential operators.

We like to mention that the theory of pseudo-differential operators 
is a theoretical foundation for the observations made already earlier. 
Indeed, it has been observed that the inverse of a kernel matrix is
compressible, compare for example \cite{Hack1,Hack2} and more recently
\cite{AFM} in case of $\mathcal{H}$-matrices and \cite{B96,SchW06} 
for the case of wavelet matrix compression. Especially, in the 
case of the Mat\'ern kernel, the associated pseudo-differential 
operator corresponds to a Bessel potential operator. The symbol 
of the Mat\'ern kernel is $(1+\|{\bs\xi}\|^2)^{d/2-\nu}$, which 
implies that the inverse pseudo-differential operator has the symbol 
$(1+\|{\bs\xi}\|^2)^{\nu-d/2}$. In case of $\nu = 1/2,3/2,\ldots$, 
this represents a differential operator and hence a local operator.
This explains the exponential decay \eqref{eq:exponentialDecayMatrix} 
of the inverse to the kernel matrix in this case, see also \cite{Rue}. 
Note that, the symbol is an analytic function in $r=\|{\bs\xi}\|$
for $\nu \neq1/2,3/2,\ldots$, hence the Schwartz kernel of the 
inverse kernel matrix still decreases rapidly as $r\to\infty$. We 
should finally mention that one may weaken the strong assumptions
on the kernels under consideration, if the Gevrey extension of 
the symbolic calculus for classical pseudo-differential operators 
is used as developed in \cite{KreeMonvel67,KreeKernel,Rodino92}.

\section{Efficient Computation of the Lagrange Basis}
\label{sec:LocLagrangeFunctions}
We now turn to introducing a framework to efficiently compute
the Lagrange basis. As already discussed above and 
quantified in \cref{coro:off-diagonal decay}, the kernel matrix
of the Mat\'ern kernel and more generally kernels of 
pseudo-differential operators exhibit a fast off-diagonal decay. The authors 
of \cite{hangelbroek:InverseTheoremLocalizedLagrangeFunctions,
fuselier:LocalizedLagrangeFunctions,hangelbroek:LpStabilityLagrangeBasis} 
used this fact to derive localizations of the Lagrange basis. 
First, we review the basic ideas of their approach and use them 
afterwards to derive symmetric and sparse preconditioners for the 
notoriously ill conditioned kernel matrix. 

\subsection{Cut-off and Localized Lagrange Functions}
\label{subsec:LocLagrangeFunctions}
In the case of the Mat\'ern kernel, we can exploit the exponential decay in
\eqref{eq:exponentialDecayMatrix} and either ignore coefficients 
in the expansions \eqref{eq:fullLagrangeFunction1} and 
\eqref{eq:modifiedLagrangeFunction1}, respectively, 
that are smaller than a given threshold or directly construct 
functions that satisfy the Lagrange condition only on a small neighborhood 
of $\xx_i\in X$. The resulting 
functions are called \emph{cut-off Lagrange basis} and 
\emph{localized Lagrange basis}, respectively. We now go 
into more detail.

Motivated by \eqref{eq:exponentialDecayMatrix}, we define 
the \emph{footprint} $X_i$ of the $ i $-th Lagrange basis element as 
\begin{align}\label{eq:footprint}
 X_i = \{ \xx_j \in X \; : \; \| \xx_i - \xx_j \|_2 
 \leq \kappa h_{X,\Omega} | \log( h_{X,\Omega}) | \},
\end{align}
with a free constant $ \kappa > 0 $, independent of $ i $. Then, 
the cut-off Lagrange basis element $ \chi_i^{\operatorname{co}} $ 
is given as
\begin{align*}
    \chi_i^{\operatorname{co}} (\xx) = \sum_{j \; : \; \xx_j \in X_i}
        \big[\AA^{-1}\big]_{i,j} K (\xx, \xx_j).
\end{align*}
Analogously, we have
\begin{align}\label{eq:localizedLagrangeFunction}
{\chi}_{\lambda,i}^{\operatorname{co}} (\xx) = \sum_{j \; : \; \xx_j \in X_i}
        \big[(\AA + \lambda\II)^{-1}\big]_{i,j} K (\xx, \xx_j),
\end{align}
for the modified cut-off Lagrange basis. Note that we only 
changed the number of indices entering the summation and kept 
the coefficients in the expansion the same as in the full Lagrange basis. 

The idea of the localized Lagrange basis is now to enforce 
the Lagrange condition only on the footprint, i.e., we define the
basis elements $ \chi_i^{\operatorname{loc}} $ and 
${\chi}_{\lambda,i}^{\operatorname{loc}} $ as
\begin{align*}
\chi_i^{\operatorname{loc}}\isdef 
\sum_{j \; : \; \xx_j \in X_i} \alpha_j^{(i)} K(\cdot, \xx_j), \quad \text{and}  \quad 
{\chi}_{\lambda,i}^{\operatorname{loc}} \isdef\sum_{j \; : \; \xx_j \in X_i} \beta_j^{(i)} K(\cdot, \xx_j)
\end{align*}
where the coefficient vectors $ \aalpha^{(i)}, \bbeta^{(i)} 
\in \R^{|X_i|} $ are the respective unique solutions of the linear systems
\begin{align}\label{eq:LGS1}
    \AA|_{X_i} \aalpha^{(i)} = \ee_m \quad \text{and} 
    \quad (\AA|_{X_i} + \lambda \II) \bbeta^{(i)} = \ee_m
\end{align}
with $ \AA|_{X_i} = [K(\xx_j, \xx_{\ell})]_{\xx_j, \xx_{\ell} \in X_i} $ 
being the restriction of the kernel matrix $\AA $ to indices in 
$ X_i $ and $ \ee_m \in \R^{|X_i|} $ being the unit vector whose 
$ 1 $-entry is in the (local) coordinate $ m \in \{1,\dots, |X_i|\} $ 
that relates to the (global) coordinate $ i \in \{1, \dots, N \}$.

While the cut-off Lagrange basis still requires to assemble 
the full kernel matrix we can achieve faster evaluation of the 
approximation $ s^{*,\operatorname{co}}_{\lambda} $. The localized Lagrange 
basis, on the other hand, reduces the computational cost since 
we have to solve for every $i=1,\ldots,N$ a linear system 
of the size $ |X_i| \times |X_i| $. Hence, if the number of 
points in the footprints is sufficiently small, this is cheaper than 
the cost to compute the full or cut-off Lagrange basis.

Following the representations of the interpolant or penalized least-squares 
approximation in \eqref{eq:representationInterpolantLagrangeFunctions} and
\eqref{eq:representationApproximationLagrangeFunctions}, respectively, we 
can use the cut-off and localized Lagrange functions to derive representations 
of approximations
\begin{align}\label{eq:quasiInterpolationCutOffLagrange}
    s_f^{\operatorname{co}} = \sum_{i=1}^{N} f(\xx_i) \chi^{co}_i \quad 
    \text{and} \quad s_{\lambda}^{*,\operatorname{co}} = \sum_{i=1}^{N} f_i {\chi}^{\operatorname{co}}_{\lambda,j}
\end{align}
or 
\begin{align}\label{eq:quasiInterpolationLocalizedLagrange}
    s_f^{\operatorname{loc}} = \sum_{i=1}^{N} f(\xx_i) \chi^{\operatorname{loc}}_i \quad 
    \text{and} \quad s_{\lambda}^{*,\operatorname{loc}} = \sum_{i=1}^{N} f_i{\chi}^{\operatorname{loc}}_{\lambda,i}.
\end{align}

For the localized Lagrange basis, the following estimate is known, see 
\cite{hangelbroek:AnInverseTheoremOnBoundedDomains}. Clearly, this translates 
immediately into an estimate for the difference $ s_f - s_f^{\operatorname{loc}} $.

\begin{proposition}
With the notation and assumptions of \cref{coro:off-diagonal decay}, where $ m := \nu - d/2 \in \N $, let 
the Lagrange basis element $ \chi_i $ be defined as in \eqref{eq:fullLagrangeFunction1} 
and the localized Lagrange basis element $ \chi_i^{\operatorname{loc}} $ be defined as in 
\eqref{eq:localizedLagrangeFunction} with footprint \eqref{eq:footprint}. 
Then there exists a constant $ C > 0 $ such that the error estimate
\begin{align*}
    \big\| \chi_i - \chi_i^{\operatorname{loc}}\big\|_{W^{\sigma}_p(\R^d)} 
    \leq C h_{X,\Omega}^{\kappa \frac{\eta}{2} +d -2m}, \quad i=1,\ldots,N,
\end{align*}
holds for $ 1 \leq p \leq \infty $ and $ \sigma < 2m - d + \frac{d}{p} $.
\end{proposition}

We now want to discuss the cost of computing the localized Lagrange basis.
Since we assume that the set of data sites $ X $ is quasi-uniform, we can bound the 
number of points in the footprint $ X_i $, independently of $ i $, by
\begin{align}\label{eq:EstimateFootprint}
|X_i| \leq C \left(\frac{\kappa}{d} \log N \right)^d
\end{align}
with a constant $ C > 0 $, independent of $ N $, if $ \frac{\kappa}{d} \log N  > 1 $.

Hence, we have to solve $ N $ dense linear systems of size 
$ \cO \big((\log^dN) \times (\log^dN) \big) $. This means, a 
naive ansatz, i.e., not using any further properties the local 
matrices exhibit, computes the coefficients for the localized 
Lagrange basis with cost $ \cO ( N \log^{3d}N) $. However, 
we emphasize that the systems are all independent of each other,
which makes it possible to efficiently parallelize the computation. Further, 
we note that we only need to keep $ \cO(N \log^dN) $ coefficients 
in memory, compared to the $ \cO(N^2) $ coefficients for 
the full Lagrange basis. This means, from a different point 
of view, we find an elegant compression of the inverse of the 
kernel matrix $ \AA $. 

\subsection{Symmetrized Version for Preconditioning}
\label{subsec:LocLagrangeFunctions2}
A drawback of the construction in
\cref{subsec:LocLagrangeFunctions} is that it is 
non-symmetric in the following sense. Collecting the 
column vectors $\aalpha^{(i)}$ and $\bbeta^{(i)}$ for all 
$i=1,\ldots,N$ in a global matrix $\BB$ yields a (sparse) matrix 
that satisfies $\AA \BB \approx \II$ or $(\AA+\lambda\II) 
\BB \approx \II$, respectively. Nonetheless, there holds 
$\BB \not=\BB^{\transpose}$ although we have $\AA =
\AA^{\transpose}$. This rules out the use of $\BB$ as a
preconditioner in the conjugate gradient method \cite{CG}. Indeed, in 
\cite{fuselier:LocalizedLagrangeFunctions} for example, 
only the
generalized minimal residual method \cite{GMRES} has been used which, however, 
is much more expensive compared to the conjugate gradient method.

We shall modify the construction in \eqref{eq:LGS1}. 
The matrix $\AA$ is symmetric and positive definite, 
so is the footprint matrix $\AA|_{X_i}$ is for any 
subset $X_i\subset X$. Therefore, the first possibility 
we have is to use the matrix square root $\AA|_{X_i}^{1/2}$ 
of the footprint matrix. Then, we solve in complete analogy 
to \eqref{eq:LGS1} the systems
\begin{align}\label{eq:LGS2}
   \AA|_{X_i}^{1/2}\widetilde\aalpha^{(i)} = \ee_m 
\end{align}
for all $i=1,\ldots,N$. Likewise to above, we then form 
a (sparse) matrix $\CC$ from all of the column vectors 
$\tilde\aalpha^{(i)}$, $1\le j\le N$, which satisfies 
$\AA^{1/2}\CC\approx\II$ and hence $\CC^{\transpose}\AA^{1/2}
\approx\II$. Thus, we also conclude $\CC^{\transpose}\AA\CC
\approx\II$ and $\CC\CC^{\transpose}$ defines a symmetric 
and positive definite preconditioner for $\AA$. 

Alternatively, we can also use the Cholesky decomposition 
$\AA|_{X_i} = \LL_i \LL_i^{\transpose}$ of size $ |X_i|\times 
|X_i| $ of the kernel matrix restricted to the footprints, solving the systems
\begin{align*}
   \LL_i^{\transpose}\tilde\aalpha^{(i)} = \ee_m 
\end{align*}
for all $i=1,\ldots,N$. If we then form the (sparse) matrix $\CC$ 
from all of the column vectors $\widetilde\aalpha^{(i)}$, $i=1,\ldots,N$,
then there holds $\LL^{\transpose} \CC\approx\II$ and hence $\CC^{\transpose}
\LL\approx\II$, where $\LL\LL^{\transpose} = \AA$ is the Cholesky 
decomposition of $\AA$. Consequently, $\CC^{\transpose}\LL\LL^{\transpose}
\CC\approx\II$ and $\CC\CC^{\transpose}$ defines likewise a symmetric 
and positive definite preconditioner for $\AA$. 

In case of the original footprint systems \eqref{eq:LGS1} 
and in case of the matrix square root \eqref{eq:LGS2}, the 
ordering of the local indices $X_i$ does not affect the so 
constructed global matrices $\BB$ and $\CC$, respectively. The 
reason for this is that the solution vectors $\aalpha^{(i)}$
and $\tilde\aalpha^{(i)}$, respectively, are then just permuted 
versions of the same vectors. This is completely different in 
case of the Cholesky decomposition. Only in the case of an 
ordering of the local indices which respects the ordering 
of the global ones, $\CC$ is an upper triangular matrix
since each vector $\widetilde\aalpha^{(i)}$ has zero entries 
for all indices being larger than the respective local 
index $m$. In this case, we have hence constructed a 
localized \emph{Newton basis}, compare~\cite{Pazouki2011}.

In complete analogy, we can also address the regularized
kernel matrices
\begin{align*}
(\AA+\lambda\II)|_{X_i}^{1/2}\tilde\bbeta^{(i)} = \ee_m
\quad\text{or}\quad
\LL_j^{\transpose}\tilde\bbeta^{(i)} = \ee_m,
\end{align*}
where $\LL_i\LL_i^{\transpose}$ is the Cholesky decomposition 
of the regularized kernel matrix $\AA|_{X_i}+\lambda\II$. 
This yields a symmetric and positive definite preconditioner 
for $\AA+\lambda\II$. We should finally mention that the 
computation of the matrix $\CC\CC^{\transpose}$ should be 
avoided, meaning that the matrix-vector multiplication 
$\zz=\CC\CC^{\transpose}\xx$ is computed in two steps in 
accordance with $\yy = \CC^{\transpose}\xx$ and $\zz = \CC\yy$.

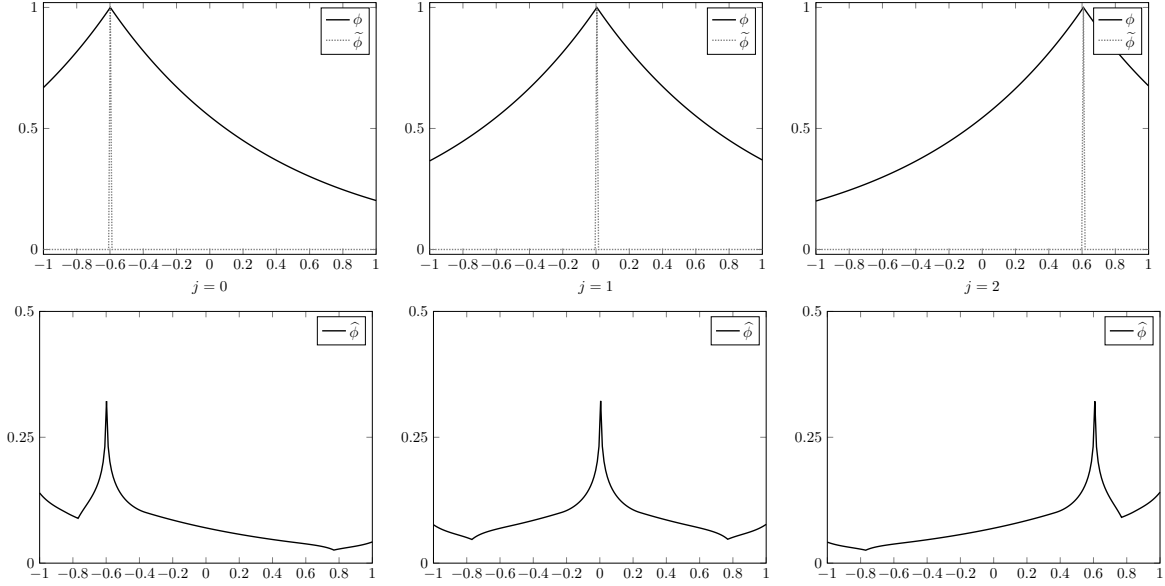
\begin{figure}[htb]
\begin{center}
\scalebox{0.52}{
\begin{tikzpicture}
\begin{axis}[width=10cm, height=8cm, xmin = -1, xmax=1,
 ymin=-0.02, ymax=1.02, ylabel={}, xlabel ={$j=0$}, ytick={0,0.5,1}]
\addplot[color=black, line width=1pt]table[x index={0},y index = {41}]{%
./Images/PK.txt};
\addlegendentry{$\phi$};
\addplot[color=gray, densely dotted, line width=1pt]
  table[x index={0},y index = {41}]{%
./Images/PsIK.txt};
\addlegendentry{$\widetilde{\phi}$};
\end{axis}
\end{tikzpicture}}
\scalebox{0.52}{
\begin{tikzpicture}
\begin{axis}[width=10cm, height=8cm, xmin = -1, xmax=1,
 ymin=-0.02, ymax=1.02, ylabel={}, xlabel ={$j=1$},ytick={0,0.5,1}]
\addplot[color=black, line width=1pt]table[x index={0},y index = {101}]{%
./Images/PK.txt};
\addlegendentry{$\phi$};
\addplot[color=gray, densely dotted, line width=1pt]
  table[x index={0},y index = {101}]{%
./Images/PsIK.txt};
\addlegendentry{$\widetilde{\phi}$};
\end{axis}
\end{tikzpicture}}
\scalebox{0.52}{
\begin{tikzpicture}
\begin{axis}[width=10cm, height=8cm, xmin = -1, xmax=1,
 ymin=-0.02, ymax=1.02, ylabel={}, xlabel ={$j=2$},ytick={0,0.5,1}]
\addplot[color=black, line width=1pt]table[x index={0},y index = {161}]{%
./Images/PK.txt};
\addlegendentry{$\phi$};
\addplot[color=gray, densely dotted, line width=1pt]
  table[x index={0},y index = {161}]{%
./Images/PsIK.txt};
\addlegendentry{$\widetilde{\phi}$};
\end{axis}
\end{tikzpicture}}
\scalebox{0.52}{
\begin{tikzpicture}
\begin{axis}[width=10cm, height=8cm, xmin = -1, xmax=1,
 ymin=0, ymax=0.5, ylabel={}, ytick={0,0.25,0.5}]
\addplot[color=black, line width=1pt]table[x index={0},y index = {41}]{%
./Images/PsIsqrtK.txt};
\addlegendentry{$\widehat{\phi}$};
\end{axis}
\end{tikzpicture}}
\scalebox{0.52}{
\begin{tikzpicture}
\begin{axis}[width=10cm, height=8cm, xmin = -1, xmax=1,
 ymin=0, ymax=0.5, ylabel={},ytick={0,0.25,0.5}]
\addplot[color=black, line width=1pt]table[x index={0},y index = {101}]{%
./Images/PsIsqrtK.txt};
\addlegendentry{$\widehat{\phi}$};
\end{axis}
\end{tikzpicture}}
\scalebox{0.52}{
\begin{tikzpicture}
\begin{axis}[width=10cm, height=8cm, xmin = -1, xmax=1,
 ymin=0, ymax=0.5, ylabel={},ytick={0,0.25,0.5}]
\addplot[color=black, line width=1pt]table[x index={0},y index = {161}]{%
./Images/PsIsqrtK.txt};
\addlegendentry{$\widehat{\phi}$};
\end{axis}
\end{tikzpicture}}
\caption{\label{fig:kernelBasesViz}
Basis of kernel translates, localized Lagrange basis for \(H^1(\Rbb)\) (top row).
The bottom row shows the corresponding orthogonal functions obtained with the
square root of the inverse Gramian \({\bs A}^{-1/2}\).
}
\end{center}
\end{figure}

In order to give a visual idea the of basis of kernel translates,
the corresponding Lagrange basis and the orthogonal basis obtained 
from the inverse square root of the kernel matrix \({\bs A}\),
we consider the Sobolev space 
$H^1(\mathbb{R})$, equipped with the usual norm 
\[
\|u\|_{H^1(\Rbb)}^2 
= \|u\|_{L^2(\Rbb)}^2 + \|u'\|_{L^2(\Rbb)}^2.
\] 
Its reproducing kernel is given by $K(x,y) = \frac{1}{2}\exp(|x-y|)$,
see \cite{wendland:ScatteredDataApproximation} for instance.
\cref{fig:kernelBasesViz} shows different basis elements and the
corresponding Lagrange basis elements. The bottom row shows the 
associated orthonormal basis, each for \(N=200\) equidistant 
data sites in $[-1,1]$.

\section{Dual Pairs in Samplet Coordinates}\label{sec:Samplets}
In \cref{subsec:LocLagrangeFunctions}, we leveraged the exponential decay 
of the entries of the inverse kernel matrix \eqref{eq:exponentialDecayMatrix} 
to construct dual pairs and quasi-interpolations to the given data and speed-up the 
evaluation of these. In this section, we pursue a multiresolution approach based
on a samplet basis embedded into a reproducing kernel Hilbert space. 
This way, we obtain a sparse representation of the kernel matrix
and a corresponding direct solver. We emphasize, however, 
that the resulting approximations will most likely not satisfy any Lagrange condition. 
Nevertheless, we will be able to define \emph{quasi-interpolations} 
in the style of \eqref{eq:quasiInterpolationCutOffLagrange} or 
\eqref{eq:quasiInterpolationLocalizedLagrange}. For the reader's convenience, we
start by briefly recalling the construction of a samplet basis. 

\subsection{Samplet Basis Construction}\label{subsec:SampletBasisConstruction}
The first building block in the samplet construction is a 
\emph{cluster tree} for the set \(X\). This is a tree
$\mathcal{T}$ with root \(X\) and each node \(\tau\in\Tcal\) is the
disjoint union of its children. We refer to \cite{HM2} for details.
The cluster tree \(\Tcal\) directly induces a support based hierarchical
clustering of the subspace 
\(\cH_X'\isdef\operatorname{span}\{\delta_{\xx_1},\ldots,\delta_{\xx_N}\}\subset\Hcal'\)
spanned by the Dirac-$\delta$-distributions supported at the data sites in \(X\).
With a slight abuse of notation, we will refer to this cluster tree also by \(\Tcal\)
and to its elements by \(\tau\). Based on the multilevel hierarchy generated
by the cluster tree, we shall now construct a multiresolution basis for \(\Hcal_X'\).
We start by introducing a \emph{two-scale 
transform} between basis elements on a cluster $\tau$ of level $j$. To this 
end, we represent scaling distributions $\bs{\Phi}_{j}^{\tau}=\{ \varphi_{j,k}^{\tau}\}$
and samplets ${\bs\Sigma}_{j}^{\tau}=\{\sigma_{j,k}^{\tau}\}$ recursively as 
linear combinations of the scaling distributions ${\bs \Phi}_{j+1}^{\tau}$ of 
$\tau$'s child clusters. This amounts to the \emph{refinement relations}
\[
\varphi_{j,k}^{\tau}=\sum_{\ell=1}^{n_{j+1}^\tau}q_{j,\Phi,\ell,k}^{\tau}
\varphi_{j+1,\ell}^{\tau}
\ \text{and}\ 
\sigma_{j,k}^{\tau}=\sum_{\ell=1}^{n_{j+1}^\tau}q_{j,\Sigma,\ell,k}^{\tau}
\varphi_{j+1,\ell}^{\tau}\ \text{with}\ n_{j+1}^\tau\isdef|
{\bs \Phi}_{j+1}^{\tau}|,
\]
which may be written in matrix notation as
\begin{equation}\label{eg:refinementRelation}
   [ {\bs \Phi}_{j}^{\tau}, {\bs \Sigma}_{j}^{\tau} ] 
 = {\bs \Phi}_{j+1}^{\tau}
 {\bs Q}_j^{\tau}=
 {\bs \Phi}_{j+1}^{\tau}
 \big[ {\bs Q}_{j,\Phi}^{\tau},{\bs Q}_{j,\Sigma}^{\tau}\big].
\end{equation}

To obtain vanishing moments and orthonormality, the transform 
\({\bs Q}_{j}^{\tau}\) is derived from an orthogonal decomposition of the 
\emph{moment matrix} ${\bs M}_{j+1}^{\tau}\in\Rbb^{m_q\times n_{j+1}^\tau}$,
given by
\[
  {\bs M}_{j+1}^{\tau}\isdef
  \begin{bmatrix}({\bs x}^{\bs 0},\varphi_{j+1,1})_\Omega&\cdots&
  ({\bs x}^{\bs 0},\varphi_{j+1,n_{j+1}^\tau})_\Omega\\
  \vdots & & \vdots\\
  ({\bs x}^{\bs\alpha},\varphi_{j+1,1})_\Omega&\cdots&
  ({\bs x}^{\bs\alpha},\varphi_{j+1,n_{j+1}^\tau})_\Omega
  \end{bmatrix}=
  [({\bs x}^{\bs\alpha},{\bs \Phi}_{j+1}^{\tau})_\Omega]_{|\bs\alpha|\le q}.
\]
Herein, $m_q=\binom{q+d}{d}$ denotes the dimension of the 
space \(\Pcal_q(\Omega)\) of total degree polynomials.

In the original construction by Tausch and White \cite{TW03}, the matrix 
\({\bs Q}_{j}^{\tau}\) is obtained from the singular value decomposition of 
\({\bs M}_{j+1}^{\tau}\). For the construction of samplets, we follow the idea
from \cite{AHK14} and rather employ the QR decomposition, which results in 
samplets with an increasing number of vanishing moments. We compute
\begin{equation}\label{eg:QR} 
  ({\bs M}_{j+1}^{\tau})^{\transpose}  = {\bs Q}_j^\tau{\bs R}
  \defis\big[{\bs Q}_{j,\Phi}^{\tau} ,
  {\bs Q}_{j,\Sigma}^{\tau}\big]{\bs R}
 \end{equation}
The moment matrix for the cluster's scaling distributions and samplets is now 
given by 
\begin{equation*}
  \begin{aligned}
  \big[{\bs M}_{j,\Phi}^{\tau}, {\bs M}_{j,\Sigma}^{\tau}\big]
  &= \left[({\bs x}^{\bs\alpha},[{\bs \Phi}_{j}^{\tau},
  	{\bs \Sigma}_{j}^{\tau}])_\Omega\right]_{|\bs\alpha|\le q}\\
  &= \left[({\bs x}^{\bs\alpha},{\bs \Phi}_{j+1}^{\tau}[{\bs Q}_{j,\Phi}^{\tau}
  , {\bs Q}_{j,\Sigma}^{\tau}])_\Omega
  	\right]_{|\bs\alpha|\le q}
  = {\bs M}_{j+1}^{\tau} [{\bs Q}_{j,\Phi}^{\tau},{\bs Q}_{j,\Sigma}^{\tau}]
  = {\bs R}^{\transpose}.
  \end{aligned}
\end{equation*}
Since ${\bs R}^{\transpose}$ is a lower triangular matrix, the first $k-1$ entries
in its $k$-th column are zero. This amounts to $k-1$ vanishing moments for the
$k$-th distribution generated by the orthogonal transform
${\bs Q}_{j}^{\tau}=[{\bs Q}_{j,\Phi}^{\tau} , {\bs Q}_{j,\Sigma}^{\tau} ]$. 
Defining the first $m_{q}$ distributions as scaling distributions and the 
remaining ones as samplets, we obtain samplets with vanishing moments at least
of order $q+1$. 

For leaf clusters, we define the scaling distributions by Dirac-$\delta$-distributions
supported at the points \({\bs x}_i\), i.e.,
${\bs \Phi}^{\tau}\isdef\{ \delta_{{\bs x}_i} : {\bs x}_i\in\tau,
\,\tau\in\Lcal(\Tcal) \}$.
The scaling distributions of all clusters on a specific level $j$ then generate
the spaces
\begin{equation}\label{HM_eg:Vspaces}
	\Xcal_{j}\isdef \spn\big\{ \varphi_{j,k}^{\tau} : 
 k\in I_j^{\Phi,\tau},\ \operatorname{level}(\tau)=j \big\}.
\end{equation}
In contrast, the samplets span the detail spaces
\begin{equation}\label{HM_eg:Wspaces}
\Scal_{j}\isdef
	\spn\big\{ \sigma_{j,k}^{\tau} : k\in I_j^{\Sigma,\tau},\ 
	\operatorname{level}(\tau)=j\big\}.
\end{equation}
Combining the scaling distributions of the root cluster with all clusters' 
samplets gives rise to the samplet basis
\begin{equation}\label{HM_eg:Wbasis}
  {\bs \Sigma}\isdef{\bs \Phi}^{X} 
  	\cup \bigcup_{\tau \in\Tcal} {\bs \Sigma}^{\tau}.
\end{equation}

\subsection{Samplet Representation of the Lagrange Basis}
\label{subsec:SampletApprox}
The samplet basis introduced in the previous subsection can be represented
according to
\begin{equation}\label{eq:samplets}
\sigma_{j,k}=\sum_{i=1}^N[\TT]_{(j,k),i}\delta_{{\bs x}_i},
\quad k\in\nabla_j,j=0,\ldots,J,
\end{equation}
where \(\nabla_j\) denote appropriate index sets.
Here, the matrix $\TT\in\mathbb{R}^{N\times N}$ describes the samplet 
transform. It is an orthogonal matrix and can be applied with linear 
cost by means of the fast samplets transform, see \cite{HM2,FWT}. 
The key feature is that they have vanishing moments of order $q+1$, 
meaning that
\begin{equation}\label{eq:vanishing}
 |(f,\sigma_{j,k})_\Omega|\le\sqrt{|\tau|} \bigg(\frac{d}{2}\bigg)^{q+1}
  	\frac{\diam(\tau)^{q+1}}{(q+1)!}\|f\|_{C^{q+1}(O)}
\end{equation}
provided that $f\in C^{q+1}(O)$, where \(O\subset\Omega\) is
any open set containing the samplet's support. Here, $\tau$ denotes the 
cluster where the samplet is supported and $|\tau|$ denotes the number 
of data sites contained in the cluster.

The samplet basis \eqref{eq:samplets} gives rise to a multiresolution 
basis in $\cH_X$ by means of the Riesz isometry \(\mathcal{J}\colon\Hcal'\to\Hcal\)
according to
\[
\psi_{j,k}\isdef\mathcal{J}\sigma_{j,k}= \sum_{i=1}^N [\TT]_{(j,k),i} K(\cdot,{\bs x}_i),
\quad k\in\nabla_j,\ j=0,\ldots, J.
\]
We call these functions \emph{embedded samplets}.
The associated Lagrange basis is given by
\[
\widetilde{\psi}_{j,k}=\sum_{j',k'} [\AA_\Sigma^{-1}]_{(j,k),(j',k')} \psi_{j',k'},
\quad k\in\nabla_j,\ j=0,\ldots, J.
\]
Here, the kernel matrix $\AA_\Sigma\in\mathbb{R}^{N\times N}$
has to be taken in samplet coordinates, i.e., we have
\begin{equation}\label{eq:samplet_matrix}
  \AA _\Sigma\isdef \TT\AA\TT^{\transpose}.
\end{equation}
Note that there holds
\[
\big\langle({\bs\Phi}{\bs U})^\intercal,
	\widetilde{\bs\Phi}{\bs U}
\big\rangle_{\Hcal}={\bs I}
\]
for any orthogonal matrix \({\bs U}\) and any dual pair \(\big({\bs\Phi},\widetilde{\bs\Phi}\big)\). 
In particular, since the samplet transform  ${\bf T}$ is orthogonal, 
the basis \(\{\psi_{j,k}\}\) corresponds to the samplet transformed basis
of kernel translates, while \(\big\{\widetilde{\psi}_{j,k}\big\}\)
is the samplet transformed Lagrange 
basis.

\begin{figure}[htb]
\begin{center}
\scalebox{0.52}{
\begin{tikzpicture}
\begin{axis}[width=10cm, height=8cm, xmin = -1, xmax=1,
 ymin=-1, ymax=2.5, ylabel={}, xlabel ={$j=0$}, ytick={-1,0,2.5}]
\addplot[color=black, line width=1pt]table[x index={0},y index = {3}]{%
./Images/PTK.txt};
\addlegendentry{$\psi$};
\addplot[color=gray, densely dotted, line width=1pt]
  table[x index={0},y index = {3}]{%
./Images/PT.txt};
\addlegendentry{$\widetilde{\psi}$};
\end{axis}
\end{tikzpicture}}
\scalebox{0.52}{
\begin{tikzpicture}
\begin{axis}[width=10cm, height=8cm, xmin = -1, xmax=1,
 ymin=-0.5, ymax=0.5, ylabel={}, xlabel ={$j=1$},ytick={-0.5,0,0.5}]
\addplot[color=black, line width=1pt]table[x index={0},y index = {9}]{%
./Images/PTK.txt};
\addlegendentry{$\psi$};
\addplot[color=gray, densely dotted, line width=1pt]
  table[x index={0},y index = {9}]{%
./Images/PT.txt};
\addlegendentry{$\widetilde{\psi}$};
\end{axis}
\end{tikzpicture}}
\scalebox{0.52}{
\begin{tikzpicture}
\begin{axis}[width=10cm, height=8cm, xmin = -1, xmax=1,
 ymin=-0.5, ymax=0.5, ylabel={}, xlabel ={$j=2$},ytick={-0.5,0,0.5}]
\addplot[color=black, line width=1pt]table[x index={0},y index = {15}]{%
./Images/PTK.txt};
\addlegendentry{$\psi$};
\addplot[color=gray, densely dotted, line width=1pt]
  table[x index={0},y index = {15}]{%
./Images/PT.txt};
\addlegendentry{$\widetilde{\psi}$};
\end{axis}
\end{tikzpicture}}
\scalebox{0.52}{
\begin{tikzpicture}
\begin{axis}[width=10cm, height=8cm, xmin = -1, xmax=1,
 ymin=-1, ymax=2.5, ylabel={}, xlabel ={$j=0$}, ytick={-1,0,2.5}]
\addplot[color=black, line width=1pt]table[x index={0},y index = {3}]{%
./Images/PTR.txt};
\addlegendentry{$\widehat{\psi}$};
\end{axis}
\end{tikzpicture}}
\scalebox{0.52}{
\begin{tikzpicture}
\begin{axis}[width=10cm, height=8cm, xmin = -1, xmax=1,
 ymin=-0.5, ymax=0.5, ylabel={}, xlabel ={$j=1$},ytick={-0.5,0,0.5}]
\addplot[color=black, line width=1pt]table[x index={0},y index = {9}]{%
./Images/PTR.txt};
\addlegendentry{$\widehat{\psi}$};
\end{axis}
\end{tikzpicture}}
\scalebox{0.52}{
\begin{tikzpicture}
\begin{axis}[width=10cm, height=8cm, xmin = -1, xmax=1,
 ymin=-0.5, ymax=0.5, ylabel={}, xlabel ={$j=2$},ytick={-0.5,0,0.5}]
\addplot[color=black, line width=1pt]table[x index={0},y index = {15}]{%
./Images/PTR.txt};
\addlegendentry{$\widehat{\psi}$};
\end{axis}
\end{tikzpicture}}
\caption{\label{fig:sampletViz}
\(H^1(\Rbb)\)-embedded primal and 
dual scaling distribution (top left), samplet on level \(j=1\) (top middle)
and samplet on level \(j=2\) (top right) for \(N=200\) equidistant 
data sites and \(q+1=3\) vanishing moments.
The bottom row shows the corresponding orthogonal samplets
obtained with the square root of the inverse kernel matrix \({\bs A}_\Sigma^{-1/2}\).
}
\end{center}
\end{figure}
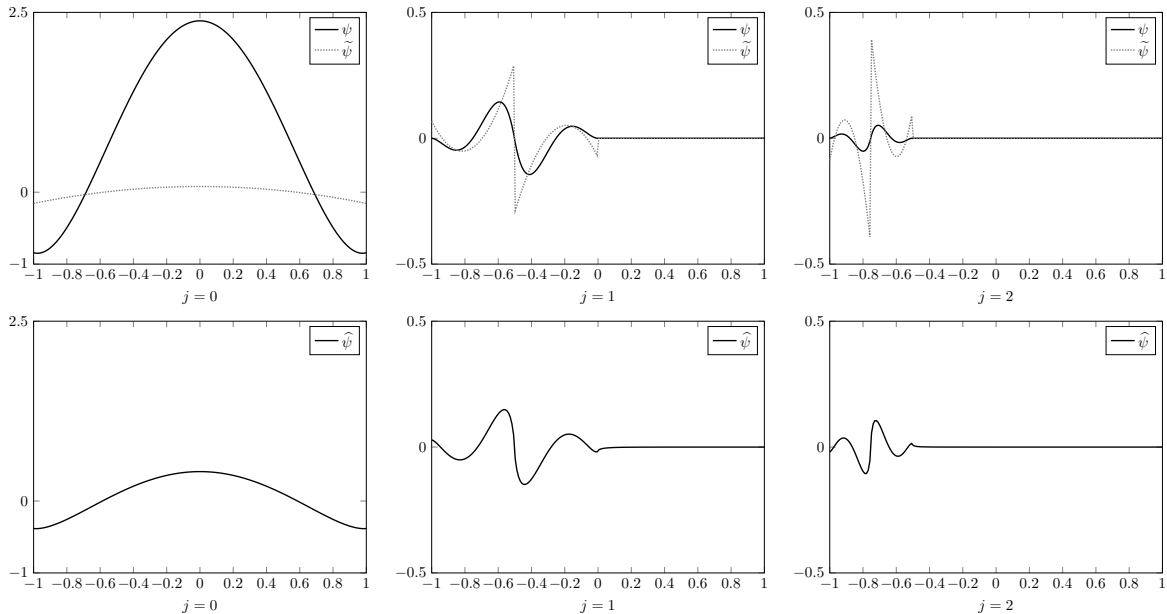

Similar to the single scale bases, we want to give
a visual idea of the {embedded samplet basis}, 
their dual basis and the corresponding orthogonal basis.
We again consider the Sobolev space 
$H^1(\mathbb{R})$. The top row of \cref{fig:sampletViz} shows an embedded 
scaling distribution (left 
plot) and two embedded samplets (middle and right plots) with \(q+1=3\) 
vanishing moments, constructed for \(N=200\) equidistant 
data sites. The bottom row shows the corresponding elements
of the orthogonal basis obtained from the inverse square root
of the inverse kernel matrix \({\bs A}_{\Sigma}^{-1/2}\). Different from the single scale
case, the orthogonal basis remains a local basis in samplet coordinates,
cp.\ \cref{fig:kernelBasesViz}.
\subsection{Kernel Matrix Compression and Associated Algebra}
\label{subsec:compression}
\begin{figure}[hbt]
\begin{center}
\includegraphics[width=0.9\textwidth]{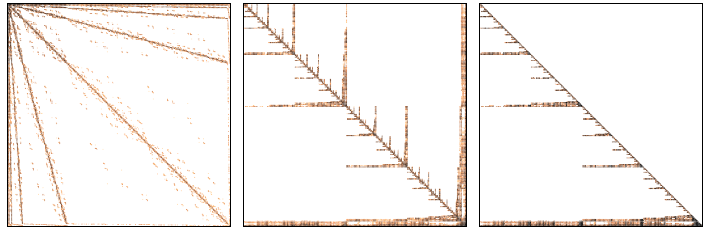}
\caption{\label{fig:S-compressed}
Typical matrix pattern of the samplet compressed kernel matrix for \(d=2\) 
(left), its reordering by means of nested dissection (middle), and the 
associated Cholesky factor (right).}
\end{center}
\end{figure}

Based on the asymptotic smoothness \eqref{eq:decay} of the 
kernel function under consideration, the kernel matrix $\AA$ 
from \eqref{eq:samplet_matrix} can be compressed in samplet 
coordinates, such that only $\mathcal{O}(N\log N)$
matrix entries remain, see \cite{HM2}. The resulting pattern 
is found in the left plot of \cref{fig:S-compressed} for uniform
random points on the unit square. 
According to \cref{subsec:pseudos}, we know that the inverse 
kernel matrix is also compressible in samplet coordinates, using 
the same matrix pattern. Therefore, it remains to compute the 
inverse in an efficient way. One approach would be the use the 
Newton-Schulz iteration to compute an approximate matrix inverse, 
see \cite{householder,schulz} for details. However, it turns out 
that selected inversion is much more efficient in terms of 
computing times, compare \cite{HMSS}. This issues from the fact
that the samplet compressed kernel matrix can be reordered by means 
of nested selection \cite{Geo73,LRT79} such that the associated 
Cholesky factors have low fill-in. The reordered kernel matrix 
is found in the middle plot of \cref{fig:S-compressed}
while the associated Cholesky factor is found in the right plot. In
the plots, lighter blocks correspond to less entries.
Selected inversion computes the inverse then based on this rather 
sparse Cholesky factor, see \cite{selinv}. The application of 
the samplet compressed, inverse kernel matrix to a data vector 
$\ff$ is computable in complexity $\mathcal{O}(N\log N)$, 
independent of the data vector being given in the samplet basis 
or -- in view of the fast samplet transform -- as point evaluations.

We finally like to mention that the rapid computation of the matrix 
power ${\bs A}^\gamma$ for $\gamma\in\mathbb{R}$ is also possible 
by means of its contour integral representation \cite{Trefethen}. 
We refer the reader to \cite{HMSS} for all the details and 
numerical results.

\section{Numerical Tests}\label{sec:Numerics}
We now turn to testing the methods introduced above. Although 
the motivation to study the construction of the localized Lagrange 
and the samplet basis was to obtain dual pairs that allow faster 
evaluation of the quasi-interpolants, both methods can be interpreted 
as compression schemes of the (inverse) kernel matrix. The compression 
power then can be directly translated into the approximation error of the 
interpolant by the respective quasi-interpolants. Hence, one part of 
the numerical tests we want to discuss is how good the matrices 
obtained from either compression method approximate the real, 
unobtainable inverse of the kernel matrix. The second part is 
dedicated to showcasing the power of the novel symmetric 
preconditioner introduced in \cref{subsec:LocLagrangeFunctions2}.
This is especially done by an application relevant problem in 
three spatial dimensions.

\subsection{Setup}
Since 
the theory of localized Lagrange functions was mainly developed for 
Mat\'ern kernels $ K = \Phi_{\nu} $, we will mostly focus on these 
kernel functions, albeit with varying parameter $ \nu = 0.5, 1.0, 1.5 $. 
We again refer to \cref{tab:exampleMatern} for explicit formulas, however, 
for the numerical tests, we choose the length scale parameter 
$ \delta = 0.1$, i.e., we consider the radial function $ \Phi_{\nu}(\cdot / \delta) $. 
This could also be considered 
for the theoretical results presented above, but only leading to 
different, $ \delta $-dependent constants.

The set of data sites $X$ is drawn randomly from a uniform distribution
on $[0,1]^2$ and we select $ N = 1000, 10\,000, 
100\,000 $ using the random number generator 
from the \texttt{C++} standard library. The fill distance $ h_{X,\Omega}$, 
that is important for determining the footprint \eqref{eq:footprint}, is then 
estimated by using a clustering algorithm. Although most of the theoretical 
results assume quasi-uniform sets of sites, the numerical examples we
conducted show that in practice this restriction can be dropped provided
that a small regularization term is added. In our experiments, we consider 
the regularized systems $ \MM + \lambda \II $ with $ \lambda = 10^{-6}N$.
Indeed, it is customary to regularize kernel systems since even for a low 
number of quasi-uniform points the unregularized kernel matrices are 
notoriously ill conditioned.

The samplet matrix compression is implemented as described in \cite{HM2}.
We compute the samplet compressed matrix using \(q+1=6\) vanishing moments.
For the fast assembly of the samplet compressed matrix, we exploit the 
fast multipole method along as proposed in \cite{HM2}. The polynomial degree 
for the degenerate kernel expansion is set to 10 while the cut-off parameter 
is chosen as \(1/2\) (see \cite{HM2} for details). After computing the 
compressed matrix, we apply different sizes of the a-posteriori 
compression parameter to steer the compression rate of the samplet 
compressed matrix. Then, we apply a sparse Cholesky decomposition 
using \texttt{METIS} \cite{METIS}.

All computations have been performed at the Centro Svizzero di Calcolo 
Scientifico (CSCS) on up to 20 nodes of the Alps cluster with two 
AMD EPYC 7742 @2.25GHz CPUs and up to 512GB of main memory,
resulting in a usage of up to $2560$ cores.

\subsection{Comparing the Compression}
We start by providing results on the compression of the inverse matrix. 
We measure the error in the following way. As explained in 
\cref{subsec:LocLagrangeFunctions2}, we collect the column 
vectors $ \bbeta^{(i)} $ of the matrix $ ( \AA + \lambda \II)^{-1} $ 
in a global (sparse) matrix $ \BB $. Since this should be a good 
approximation to the inverse of $ (\AA + \lambda \II) $, the 
error measure
\begin{align*}
   \operatorname{error} = {\| (\AA + \lambda \II) \BB - \II \|_2},
\end{align*}
would be best suited to compare the methods. Note that the spectral
error is a relative error due to \(\|\II\|_2=1\).
To compute the spectral error we use 200 power iterations.

In the theoretical results from 
\cref{subsec:LocLagrangeFunctions}, the localization of 
the Lagrange basis leads to a matrix compression and 
the resulting error depends on the parameter $ \kappa $ 
which influences the size of the footprint. However, to 
compare the numerical results with the samplet compression,
we decide to use a different parameter. Since we are mainly 
interested in how good compressed matrices for the single 
methods approximate the inverse, we consider the \emph{compression rate} 
\begin{align*}
    \operatorname{compression\ rate} = \frac{\text{number of non-zero entries}}{N^2}.
\end{align*}
This seems to be the right reference value to allow a fair 
comparison between the methods, since the compression rate 
directly reflects the memory requirements for the approximate dual basis.

\begin{figure}[htb]
\begin{center}
\pgfplotstableread{
comprate nu05 nu10 nu15
0.017884          2.59537     3.02188       4.48555    
0.06416          0.103989     0.0364512     0.162932  
0.133476        0.0164087     0.00181752    0.00662425 
0.218442       0.00360557    0.000403435    0.000359195 
0.314844       0.00089673    0.000161731    4.0232e-05 
0.418368      0.000239306    2.50578e-05    1.13117e-05
0.523304      7.08501e-05    5.40728e-06    1.24218e-06
0.62532       2.35292e-05    1.88325e-06    3.01652e-07
}\loadedtableOne 

\pgfplotstableread{
comprate nu05 nu10 nu15
0.00388614    11.1651           2.68161      13.265
0.0147652     0.97808          0.259241     3.69366
0.032132      0.277261        0.0194949    0.322787
0.055357      0.101069         0.007955   0.0174833
0.0837093     0.0423513      0.00346926  0.00415171
0.11663       0.0192898      0.00160194  0.00132455
0.153457      0.00931418     0.000759778 0.000415133
0.193469      0.00467335     0.000361902 0.000136823
}\loadedtableTwo 

\pgfplotstableread{
comprate nu05 nu10 nu15
 0.000906683  62.2871           266.423    448.263
0.00354691     6.03453          15.3282    21.1785
0.00785443     2.13354          0.57323    23.9281
0.013747       0.960075       0.0877097    4.50762
0.0211548      0.497623       0.0382435   0.527353
0.0300006      0.281387       0.0207815    0.26966
}\loadedtableThree 

\begin{tikzpicture}[scale=1]
\begin{semilogyaxis}[width=0.6\textwidth, height = 0.5\textwidth, xlabel={compression rate}, ylabel={error}, legend pos= north east,
xmin = -0.05, xmax = 0.75, ymin = 0.0000001, ymax = 3000, xtick={0,0.1,0.2,0.3,0.4,0.5,0.6,0.7},
legend style={nodes={scale=0.5, transform shape}},grid]
   \addplot[mark=*,color=red] table[header=false, x=comprate,y=nu05]{\loadedtableOne};
   \addlegendentry{Footprints, $N=1000$, $\nu = 0.5$};
   \addplot[mark=square*,color=red,style=dashed] table[header=false, x=comprate,y=nu05]{\loadedtableTwo};
   \addlegendentry{Footprints, $N=10000$, $\nu = 0.5$};
   \addplot[mark=triangle*,color=red,style=dotted] table[header=false, x=comprate,y=nu05]{\loadedtableThree};
   \addlegendentry{Footprints, $N=100000$, $\nu = 0.5$};   
   \addplot[mark=*,color=blue] table[header=false,x=comprate,y=nu10] {\loadedtableOne}; 
   \addlegendentry{Footprints, $N=1000$, $\nu = 1.0$};
   \addplot[mark=square*,color=blue,style=dashed] table[header=false,x=comprate,y=nu10] {\loadedtableTwo}; 
   \addlegendentry{Footprints, $N=10000$, $\nu = 1.0$};
    \addplot[mark=triangle*,color=blue,style=dotted] table[header=false,x=comprate,y=nu10] {\loadedtableThree}; 
    \addlegendentry{Footprints, $ N=100000$, $\nu = 1.0$};     
   \addplot[mark=*,color=brown] table[header=false, x=comprate,y=nu15] {\loadedtableOne}; 
   \addlegendentry{Footprints, $N=1000$, $\nu = 1.5$};
      \addplot[mark=square*,color=brown,style=dashed] table[header=false, x=comprate,y=nu15] {\loadedtableTwo}; 
      \addlegendentry{Footprints, $N=10000$, $\nu = 1.5$};
    \addplot[mark=triangle*,color=brown,style=dotted] table[header=false, x=comprate,y=nu15] {\loadedtableThree}; 
    \addlegendentry{Footprints, $N=100000$, $\nu = 1.5$};
 \end{semilogyaxis} 
\end{tikzpicture}
\caption{\label{fig:CompErr2DFoot}
Localized Lagrange approach: Errors as a function of the compression 
rate in case of varying smoothness parameters $\nu = 0.5, 
1.0, 1.5$ and sample sizes $N = 10^3, 10^4, 10^5$.}
\end{center}
\end{figure}
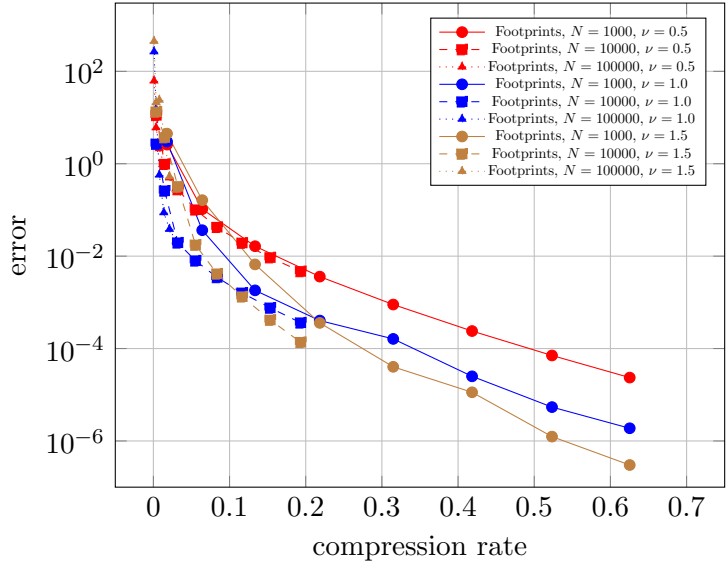

In \cref{fig:CompErr2DFoot}, we plot the error as a function 
of the compression rate for the local Lagrange approach. It is 
interesting to see that the error behaves for every choice of $ \nu $ 
qualitatively the same and the method seems to converge, achieving 
good approximation errors even with high compression. The
results indicate that the localization of 
the Lagrange basis yields a simple to implement approach, 
which gives accuracies dependent just on the compression 
rate. Indeed, looking at \cref{fig:CompErr2DFoot}, the error 
seems to be basically independent of the sample size $N$ and 
the smoothness parameter $\nu$.

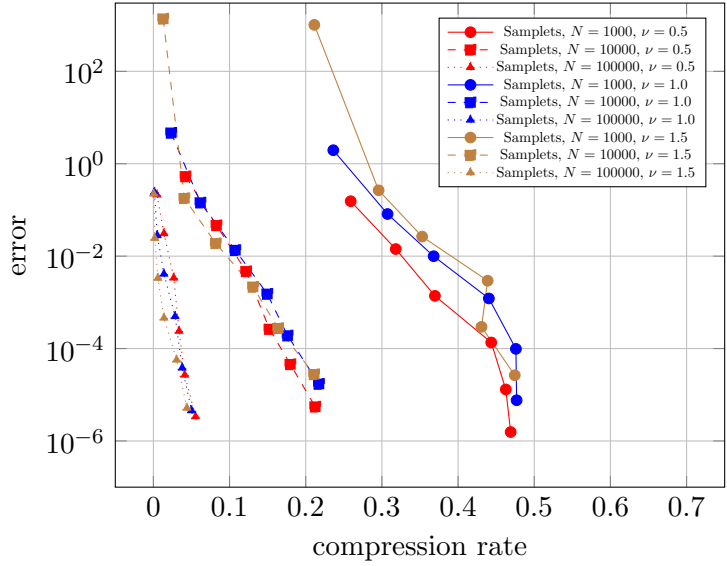
\begin{figure}[htb]
\begin{center}
\pgfplotstableread{
comprate nu05
0.259028 0.154619
0.318177 0.0142479
0.369593 0.00138242
0.443541 0.000135356
0.462635 1.30234e-05
0.46905  1.56149e-06
}\loadedtableOneA 

\pgfplotstableread{
comprate nu10
0.236208  1.96338
0.307427  0.0815755
0.367939  0.00996799
0.440417  0.00121495
0.47593   9.83953e-05
0.476869  7.60749e-06
}\loadedtableOneB 

\pgfplotstableread{
comprate nu15
0.211281 1013.54
0.295701 0.266118
0.352833 0.0265146
0.438838 0.00293831
0.430784 0.000295043
0.474394 2.64471e-05
}\loadedtableOneC 

\pgfplotstableread{
comprate nu05
0.0422183 0.532022
0.0828362 0.0464166
0.121808 0.00469084
0.152349  0.000264583
0.179837  4.56313e-05
0.212455  5.51741e-06
}\loadedtableTwoA 

\pgfplotstableread{
comprate nu10
0.0233716  4.66276
0.0617961  0.143996
0.107561   0.0134809
0.149554   0.00152196
0.17619    0.000191026
0.21732    1.71418e-05
}\loadedtableTwoB 

\pgfplotstableread{
comprate nu15
0.0128979 1374.12
0.0406888 0.179051
0.0818776 0.0188847
0.130862 0.00216178
0.163849 0.000275348
0.211032 2.76979e-05
}\loadedtableTwoC 

\pgfplotstableread{
comprate nu05
0.00470284 0.212432
0.0138995  0.0309023
0.0267956  0.00338441
0.0337277  0.000238163
0.0411822  2.68456e-05
0.0553809  3.36456e-06
}\loadedtableThreeA 

\pgfplotstableread{
comprate nu10
0.00138274 0.241469
0.00510587 0.0280089
0.0138972  0.00410215
0.0284993  0.000495636
0.0378475  3.82059e-05
0.0501932  4.53964e-06
}\loadedtableThreeB 

\pgfplotstableread{
comprate nu15
0.000557502 0.221678
0.00185265 0.0244146
0.00588102 0.00333675
0.0139051 0.000459313
0.0304777 5.67595e-05
0.0440397 5.15808e-06
}\loadedtableThreeC 

\begin{tikzpicture}[scale=1]
 \begin{semilogyaxis}[width=0.6\textwidth, height = 0.5\textwidth,xlabel={compression rate}, ylabel={error}, legend pos=north east, 
 xmin = -0.05, xmax = 0.75, ymin = 0.0000001, ymax = 3000, xtick={0,0.1,0.2,0.3,0.4,0.5,0.6,0.7},
 legend style={nodes={scale=0.5, transform shape}},grid] 
   \addplot[mark=*,color=red] table[header=false, x=comprate,y=nu05]{\loadedtableOneA};
   \addlegendentry{Samplets, $N=1000$, $\nu = 0.5$};
   \addplot[mark=square*,color=red,style=dashed] table[header=false, x=comprate,y=nu05]{\loadedtableTwoA};
   \addlegendentry{Samplets, $N=10000$, $\nu = 0.5$};
   \addplot[mark=triangle*,color=red,style=dotted] table[header=false, x=comprate,y=nu05]{\loadedtableThreeA};
   \addlegendentry{Samplets, $N=100000$, $\nu = 0.5$};   
   \addplot[mark=*,color=blue] table[header=false,x=comprate,y=nu10] {\loadedtableOneB}; 
   \addlegendentry{Samplets, $N=1000$, $\nu = 1.0$};
   \addplot[mark=square*,color=blue,style=dashed] table[header=false,x=comprate,y=nu10] {\loadedtableTwoB}; 
   \addlegendentry{Samplets, $N=10000$, $\nu = 1.0$};
    \addplot[mark=triangle*,color=blue,style=dotted] table[header=false,x=comprate,y=nu10] {\loadedtableThreeB}; 
    \addlegendentry{Samplets, $ N=100000$, $\nu = 1.0$};     
   \addplot[mark=*,color=brown] table[header=false, x=comprate,y=nu15] {\loadedtableOneC}; 
   \addlegendentry{Samplets, $N=1000$, $\nu = 1.5$};
      \addplot[mark=square*,color=brown,style=dashed] table[header=false, x=comprate,y=nu15] {\loadedtableTwoC}; 
      \addlegendentry{Samplets, $N=10000$, $\nu = 1.5$};
    \addplot[mark=triangle*,color=brown,style=dotted] table[header=false, x=comprate,y=nu15] {\loadedtableThreeC}; 
    \addlegendentry{Samplets, $N=100000$, $\nu = 1.5$};
 \end{semilogyaxis} 
\end{tikzpicture}
\caption{\label{fig:CompErr2DSmp}
Samplet matrix compression: Errors as a function of the compression 
rate in case of varying smoothness parameters $\nu = 0.5, 
1.0, 1.5$ and sample sizes $N = 10^3, 10^4, 10^5$.}
\end{center}
\end{figure}


In contrast, as seen in \cref{fig:CompErr2DSmp},
the compression rate of the samplet method improves as 
$N$ increases. The reason for this is that the samplet 
compression gives a fixed error in basically (i.e., up 
to a polylogarithmic factor) linear cost, meaning 
that the compression rate improves as $N$ increases. This 
can clearly be observed from \cref{fig:CompErr2DSmp}. Note 
that the compression rate specified in this figure is the one 
of the sparse Cholesky factor since the approximate inverse 
matrix is symmetric and only half the matrix needs to be stored.
When the whole inverse should be stored, the compression 
rate is twice the given values.

We may conclude from comparing \cref{fig:CompErr2DFoot} and
\cref{fig:CompErr2DSmp} (note that the scaling of the axes 
is the same) that samplet matrix compression is superior to
the local Lagrange basis approach, at least for large data sets. 
For small data sets, however, the local Lagrange basis approach 
approximates the inverse matrix quite well and is an 
attractive alternative since it is much easier to implement 
than the samplet matrix compression.

\subsection{Testing the Symmetric Preconditioner}
We now turn to testing the new symmetric preconditioner 
introduced in \cref{subsec:LocLagrangeFunctions2} using the 
localized Lagrange basis. We consider the same setup as before but restrict 
ourselves to a fixed amount of $N=100\,000$ data sites. 
To assess the quality of the preconditioner,
we are interested in comparing the number of (preconditioned) 
conjugate gradient
iterations (see \cite{CG} for details) when solving the linear 
system \eqref{eq:LGS} of equations for the right-hand side 
given by the vector of all ones versus the footprint parameter 
$\kappa$. To this end, the conjugate gradient method is stopped when the 
residual is of order $ 10^{-9} $.


\begin{table}[hbt]
\begin{center}
\begin{tabular}{|c|c|c|c|c|c|c|c|}\cline{3-8}
\multicolumn{2}{c}{} & \multicolumn{3}{|c|}{Cholesky decomposition} & \multicolumn{3}{c|}{Matrix square root} \\ \hline
$\kappa$ & $\# X_j$ & $\nu = 0.5$ & $\nu = 1.0$ & $\nu = 1.5$ & $\nu = 0.5$ & $\nu = 1.0$ & $\nu = 1.5$ \\ \hline
0.5 & 90 & 208 & 339 & 463 & 60 & 63 & 78 \\
1.0 & 354 & 92 & 136 & 174 & 29 & 28 & 26 \\
1.5 & 785 & 56 & 77  & 93 & 21 & 18 & 17 \\
2.0 & 1374 & 39 & 49 & 58 & 16 & 14 & 13 \\
2.5 & 2115 & 30 & 35 & 41 & 14 & 11 & 10 \\
3.0 & 3000 & 24 & 28 & 31 & 11 & 10 & 9 \\ \hline
\end{tabular}
\caption{
\label{tab:CG}
Conjugate gradient iterations in case of Cholesky decomposition and the matrix square root
versus the footprint size. The matrix square root is more stable since it
performs better for the same footprint size, but it requires also more 
computation time.}
\end{center}
\end{table}

We observe from \cref{tab:CG} that the matrix square root
of the footprints provides a better preconditioner 
compared to their Cholesky decomposition. Indeed, the number 
of iterations are by about a factor 2--3 smaller for the same
footprint size. Nonetheless, the computation of the matrix square 
root of the footprints is much more expensive compared to their 
Cholesky decomposition. First, we need only half the memory to store 
the Cholesky factor compared to the matrix square root. Second, the 
computation of the matrix square root involves the computation of 
all eigenvalues and eigenvectors by the QR-method, which is costly. 
If the problem size does not fit any more into the cache of the 
computer, we require a factor of about 100 more computation time 
per footprint for the matrix square root compared to the Cholesky 
decomposition. 

From these experiments, we may thus conclude that the Cholesky 
decomposition of the footprint matrices provides a good, symmetric 
preconditioner for large kernel matrices, being superior over 
the matrix square root approach. We will demonstrate this 
finding by a further experiment in the next subsection.

\subsection{Signed Distance Function Interpolation}

\begin{figure}[htb]
\begin{center}
\begin{tikzpicture}
\draw(-5,0)node{\includegraphics[scale=0.065,clip,trim= 1000 100 1000 100]{
  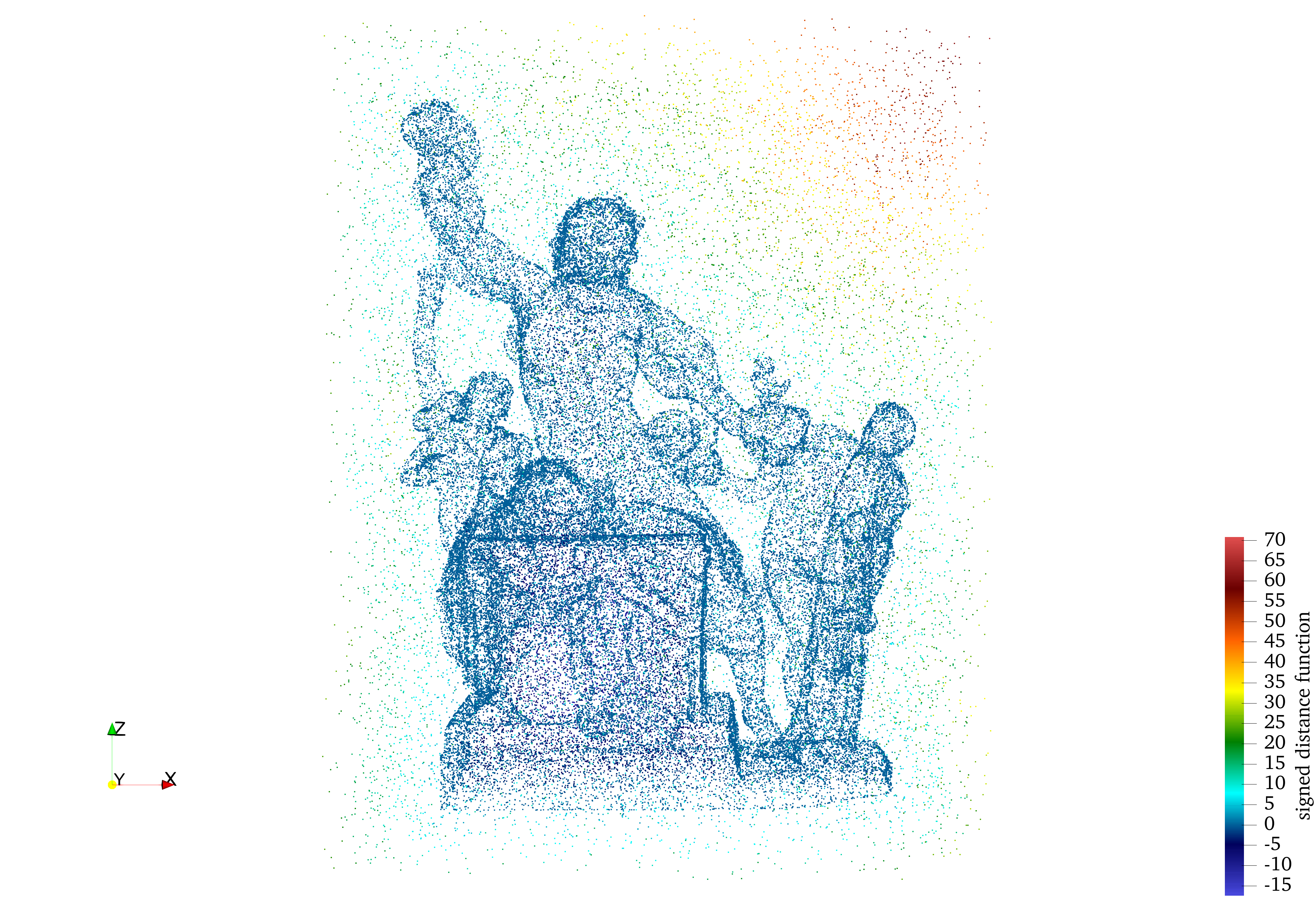}};
\draw(0,0)node{\includegraphics[scale=0.108,clip,trim= 450 90 450 90]{
  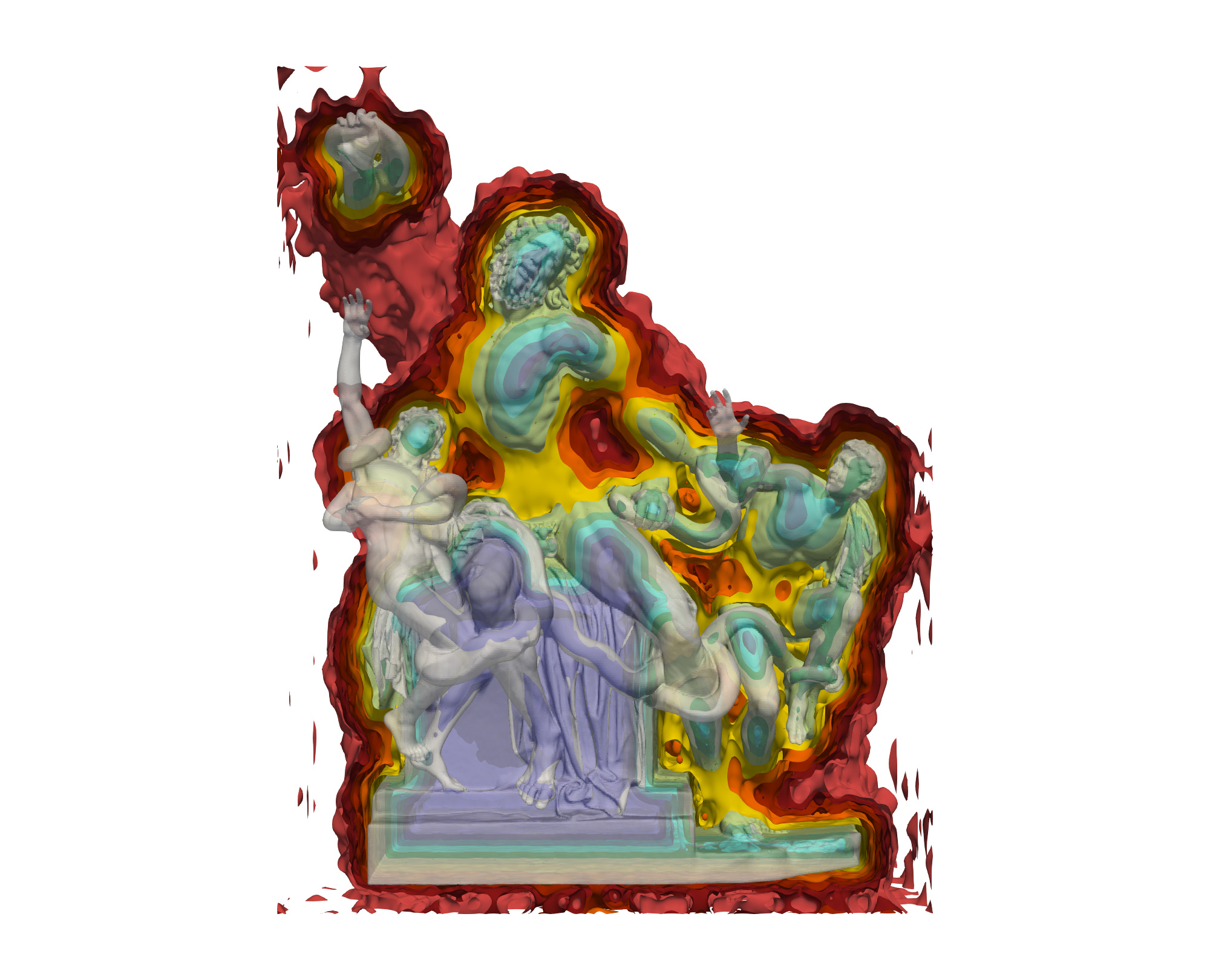}};
\draw(5,0)node{\includegraphics[scale=0.108,clip,trim= 450 90 450 90]{
  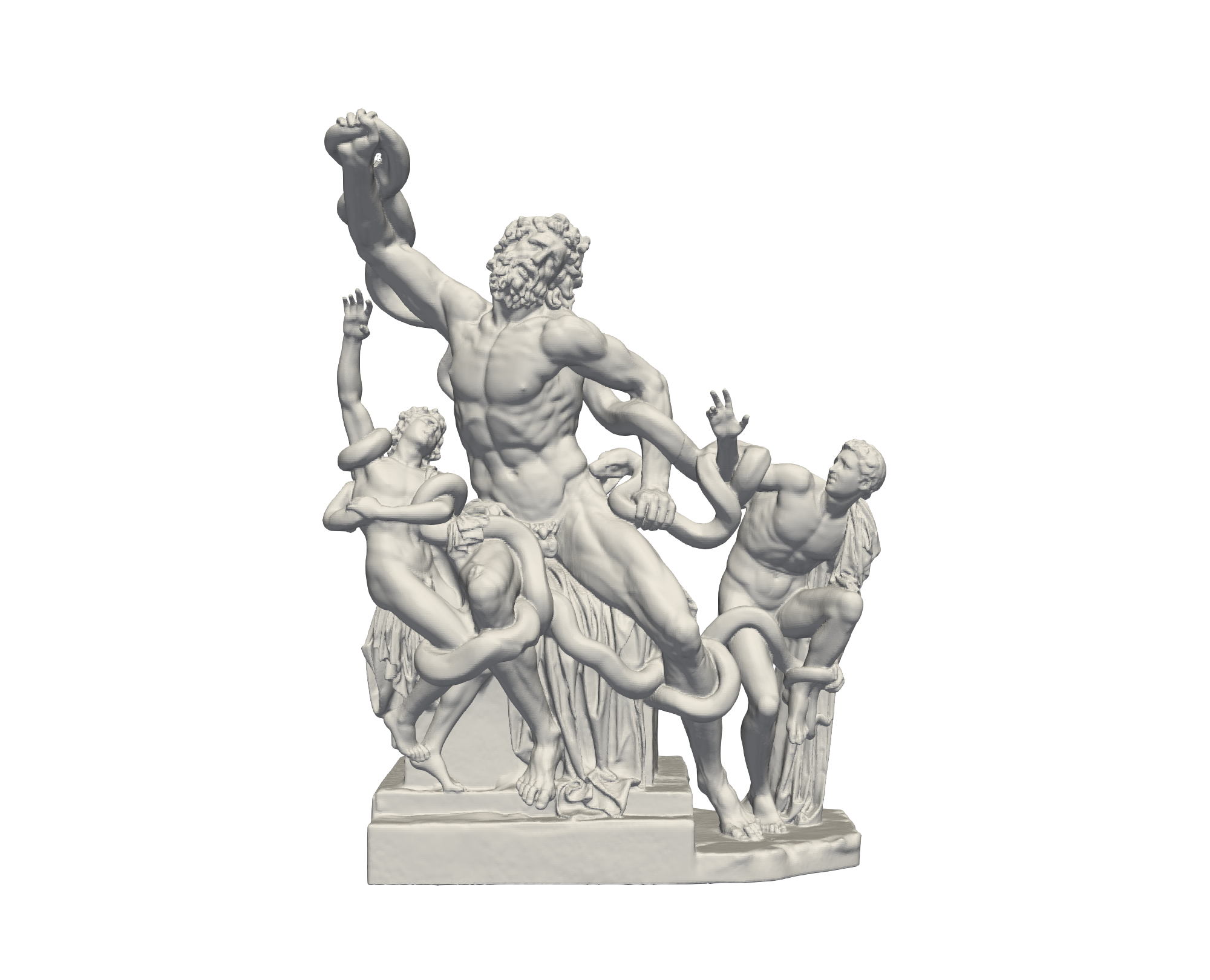}};
\end{tikzpicture}
\caption{\label{fig:SDFInterpolation}
Signed distance function interpolation
from a surface mesh of the Laokoon group: 
Subsample of the signed distance function data (left),
levelsets for the values
\(\{-5,-4,\ldots,5\}\) (middle) and zero levelset (right).}
\end{center}
\end{figure}

To demonstrate the capabilities of the developed preconditioner
also for real world applications, we consider a signed distance 
function interpolation problem in three spatial dimensions and the 
respective surface reconstruction, similarly to \cite{CarrEtal01}. 

Given a planar triangulation from a 3D scan of the Laokoon 
group (the scan is provided by the Statens Museum for Kunst), we 
generate uniform samples of the signed distance by using the about 
500\,000 vertices the surface mesh and another 250\,000 random points 
within the bounding box of the Laokoon group. This results in 
\(N=750\,000\) samples of the signed distance function in total. 
The left image in \cref{fig:SDFInterpolation} shows a uniform 
subsample of size 100\,000 of the data points. For the interpolation, 
we consider the exponential kernel \(\Phi_{1/2}\) with lengthscale 
parameter \(\delta=0.01\), where the data sites are rescaled to the 
hypercube $(0,1)^3$. The kernel matrix is compressed using a hierarchical 
matrix with adaptive cross approximation, see \cite{BR03}, to compress 
the far-field of the kernel. The achieved compression error is 
$2.37305\cdot 10^{-7}$ with a compression rate of $0.009375$. The 
kernel matrix is regularized with a regularization parameter \(\lambda=10^{-8}N\).

We use the variant of the preconditioner which is based on the Cholesky 
decomposition as it has proven to be computationally significantly cheaper 
than the one which is based on the matrix square root. We set \(\kappa=0.25\), 
which results in a compression rate of $0.000268$. The corresponding spectral 
error is \(0.657666\) and the conjugate gradient method takes only $12$ 
iterations to achieve a relative residual error of $5.33355\cdot10^{-9}$.

The interpolated signed distance function is evaluated at a uniform grid with
500 points per axis direction, resulting in 125\,000\,000 points in total.
For the evaluation, we employ the fast multipole method developed in \cite{HMQ24} 
with polynomials of total degree 4 for the degenerate kernel expansion of the
kernel in the far-field. The middle image of \cref{fig:SDFInterpolation} 
shows the levelsets for the values \(\{-5,-4,\ldots,5\}\), while the zero levelset, 
corresponding to the reconstructed surface, is found on the right hand side of 
the figure.

\section{Conclusion}\label{sec:conclusio}
Motivated by finding dual pairs for finite dimensional subspaces 
of reproducing kernel Hilbert spaces, we have discussed two methods to 
approximate the canonical dual pair. On the one hand, we used the 
exponential decay of the inverse of the Mat\'ern kernel to obtain 
the local Lagrange functions. The ideas 
are also applicable to other kernels, those that induce a 
pseudo-differential operator. 

On the other hand, we gave a short overview of the samplet 
approximation of Lagrange functions, a method that works 
also for a wider class of kernels, in particular those whose 
inverse only exhibit algebraic decay.  

Both approaches can also be used to compress the (inverse) of 
the kernel matrix. The matrix resulting from 
the localized Lagrange functions has already been used 
for preconditioning of the GMRES method, we, for the first 
time however, derived a symmetric preconditioner for the 
kernel matrix of the Mat\'ern kernel. This means that 
the preconditioner also applies to the conjugate 
gradient method.

We have presented extensive numerical tests. On the one hand, 
we have compared the compressive power of the two approaches. 
Here, we saw that only measuring the compression, 
the samplet method are clearly superior compared to the local 
Lagrange basis method. However, the latter method is very easy to 
implement and is arbitrarily scalable whereas the samplet 
method needs very sophisticated algorithms to perform. 
Finally, we have tested the new preconditioning method in 
an academic example as well as a possible real-live 
application to demonstrate its feasibility.

\bigskip\noindent
{\bf Acknowledgement.}
HH and RK were funded in parts by the Swiss National Science 
Foundation through the grant ``Adaptive Boundary Element Methods 
Using Anisotropic Wavelets'' (200021\_192041). MM was funded in 
parts by the SNSF starting grant ``Multiresolution methods for 
unstructured data'' (TMSGI2\_211684).

\appendix

\bibliographystyle{plain}
\bibliography{kempf}
\end{document}